\documentstyle[11pt,dgspp,euscript]{article}
\newtheorem{theorem}{Theorem}[section]
\newtheorem{definition}[theorem]{Definition}
\newtheorem{proposition}[theorem]{Proposition}
\newtheorem{lemma}[theorem]{Lemma}
\newtheorem{corollary}[theorem]{Corollary}
\newremark{example}[theorem]{Example}
\newremark{remark}[theorem]{Remark}
\newremark{remarks}[theorem]{Remarks}

\newcommand{\bs}{\backslash}
\newcommand{\bul}{\,\scriptscriptstyle\bullet}
\newcommand{\bve}{\bar{\varepsilon}}
\newcommand{\del}[1]{\nabla_{#1}}
\newcommand{\dg}{\dot\gamma}
\newcommand{\dps}{\displaystyle}
\newcommand{\ds}{\oplus} %direct sum 
\newcommand{\dx}{\dot{x}}
\newcommand{\ddx}{\ddot{x}}
\newcommand{\e}{\mbox{$\frak{e}$}}
\newcommand{\g}{\gamma}
\newcommand{\h}{\mbox{$\frak{h}$}}
\newcommand{\half}{\mbox{$\txs\frac{1}{2}$}}
\newcommand{\io}{\iota}
\newcommand{\Iaut}{I^{a\kern-.05em u\kern-.02em t}}
\newcommand{\Ispl}{I^{s\kern-.055em p\kern-.025em l}}
\newcommand{\lsp}{[\kern-0.15em[} %left spanning bracket 

\newcommand{\n}{\mbox{$\frak{n}$}}
\newcommand{\om}{\omega}
\newcommand{\ph}{$p\kern-.1em H\!$}
\newcommand{\quar}{\mbox{$\txs\frac{1}{4}$}}
\newcommand{\rsp}{]\kern-0.15em]} %right spanning bracket 
\newcommand{\sC}{\EuScript C}
\newcommand{\sJ}{\EuScript J}

\newcommand{\surj}{\rightarrow\kern-.82em\rightarrow}
\newcommand{\th}{\theta}
\newcommand{\tphi}{\tilde{\phi}}
\newcommand{\txs}{\textstyle}
\newcommand{\ve}{\mbox{$\varepsilon$}}
\newcommand{\vp}{\mbox{$\varphi$}}
\newcommand{\w}{\frak{w}}
\newcommand{\wO}{\widetilde{O}}
\newcommand{\z}{\frak{z}}
\newcommand{\E}{\frak{E}}
\newcommand{\F}{\mbox{$\Phi_t$}}
\newcommand{\G}{\mbox{$\Gamma$}}
\newcommand{\I}{\rsfs{I}}
\newcommand{\II}{\rsfs{I{\kern-.55em}I}}

\newcommand{\R}{\mbox{${\Bbb R}$}}
\newcommand{\U}{\frak{U}}
\newcommand{\V}{\frak{V}}
\newcommand{\Z}{\frak{Z}}

\renewcommand{\a}{\alpha}

\renewcommand{\l}{\lambda}
\renewcommand{\v}{\frak v}

\makeatletter
\newcommand{\Ad}[1]{\mathop{\operator@font Ad}\nolimits_{#1}}
\newcommand{\ad}[1]{\mathop{\operator@font ad}\nolimits_{#1}}
\newcommand{\add}[2]{\mathop{\operator@font ad}\nolimits^{\dagger}_{#1}{#2}}
\newcommand{\Aut}{\mathop{\operator@font Aut}\nolimits}
\newcommand{\diag}{\mathop{\operator@font diag}\nolimits}   
\newcommand{\End}{\mathop{\operator@font End}\nolimits}
\newcommand{\Ric}{\mathop{\operator@font Ric}\nolimits}
\newcommand{\Dspecp}{\mathop{\operator@font {\EuScript D}spec}\nolimits_\wp}
\newcommand{\specl}{\mathop{\operator@font spec}\nolimits_\ell}
\newcommand{\specp}{\mathop{\operator@font spec}\nolimits_\wp}
\newcommand{\tr}{\mathop{\operator@font tr}\nolimits}
\makeatother

\newcommand{\ns}{\normalshape} %a real hack 

\preprint{CP4}
\title{PSEUDORIEMANNIAN\\ 2-STEP NILPOTENT LIE GROUPS}
\author{Luis A. Cordero\thanks{Partially supported by Project
   XUGA20702B96, Xunta de Galicia, Spain.}}
\address{Dept. Xeometr{\'\i}a e Topolox{\'\i}a\\
        Facultade de Matem\'aticas\\
        Universidade de Santiago\\
        15706 Santiago de Compostela\\
        Spain\\
        cordero@zmat.usc.es}
\author{Phillip E. Parker\thanks{Partially supported by MEC:DGES
Program SAB1995-0757, Spain.}}
\address{Mathematics Department\\
        Wichita State University\\
        Wichita KS 67260-0033\\
        USA\\
        pparker@twsuvm.uc.twsu.edu}
\date{10 May 1999} %\draft 
\abstract{We begin a systematic study of these spaces, initially following
along the lines of Eberlein's comprehensive study of the Riemannian case.
In particular, we integrate the geodesic equation, discuss the structure of
the isometry group, and make a study of lattices and periodic geodesics.

Some major differences from the Riemannian theory appear.  There are many
flat groups ({\em versus\/} none), including Heisenberg groups.  While
still a semidirect product, the isometry group can be strictly larger than
the obvious analogue.  Everything is illustrated with explicit examples.

We introduce the notion of \ph-type, which refines Kaplan's $H\!$-type and
completes Ciatti's partial extension.  We give a general construction for
algebras of \ph-type.}
\msc{53C50}{22E25, 53B30, 53C30}

\begin{document}

\maketitle

\section{Introduction}

While there had not been much published on the geometry of nilpotent Lie
groups with a left-invariant Riemannian metric in 1990 \cite{E}, the
situation is certainly better now; see \cite{LP,Ma,Wl} and the references in
\cite{E'}.  However, there is still almost nothing extant about the more
general pseudoriemannian case.  In particular, the 2-step nilpotent groups
are nonabelian and as close as possible to being Abelian, but display a rich
variety of new and interesting geometric phenomena.  As in the Riemannian
case, one of many places where they arise naturally is as groups of
isometries acting on horospheres in certain (pseudoriemannian) symmetric
spaces.  Another is in the Iwasawa decomposition of semisimple groups with
the Killing metric, which need not be definite.  Here we begin the study of
these groups.

One motivation for our study was our observation in \cite{CP3} that there
are two nonisometric pseudoriemannian metrics on the Heisenberg group
$H_3$, one of which is flat.  This is a strong contrast to the Riemannian
case in which there is only one (up to positive homothety) and it is {\em
not\/} flat.  This is not an anomaly, as we shall see later.  We were also
inspired by the paper of Eberlein \cite{E,E'}, and followed it quite
closely in some places.  Since the published version is not identical to
the preprint, we have cited both where appropriate.

While the geometric properties of Lie groups with left-invariant definite
metric tensors have been studied extensively, the same has not occurred for
indefinite metric tensors.  For example, while the paper of Milnor \cite{Mi}
has already become a classic reference, in particular for the classification
of positive definite (Riemannian) metrics on 3-dimensional Lie groups, a
classification of the left-invariant Lorentzian metric tensors on these
groups became available only very recently \cite{CP3}.  Similarly, only a
few partial results in the line of Milnor's study of definite metrics were
previously known for indefinite metrics \cite{B,No}.  Moreover, in dimension
3 there are only two types of metric tensors:  Riemannian (definite) and
Lorentzian (indefinite).  But in higher dimensions there are many distinct
types of indefinite metrics while there is still essentially only one type
of definite metric.  This is another reason our work here has special
interest.

\medskip
By an {\em inner product\/} on a vector space $V$ we shall mean a
nondegenerate, symmetric bilinear form on $V$, generally denoted by
$\langle\,,\rangle$.  In particular, we {\em do not\/} assume that it is
positive definite.  Our convention is that $v\in V$ is timelike if
$\langle v,v\rangle > 0$, null if $\langle v,v\rangle = 0$, and spacelike
if $\langle v,v\rangle < 0$.

Throughout, $N$ will denote a connected, 2-step nilpotent Lie group with
Lie algebra \n\ having center $\z$.  (Recall that 2-step means
$[\n,\n]\subseteq\z$.) We shall use $\langle\,,\rangle$ to denote either
an inner product on \n\ or the induced left-invariant pseudoriemannian
(indefinite) metric tensor on $N$.

In Section~\ref{defex}, we give the fundamental definitions and examples
used in the rest of this paper.  The main problem encountered is that the
center $\z$ of \n\ may be degenerate:  it might contain a (totally) null
subspace.  We shall see that this possible degeneracy of the center causes
the essential differences between the Riemannian and pseudoriemannian cases.
At the end is our pseudoriemannian \ph-type, and a kind of equivalence with
$H\!$-type.  Our \ph-type actually refines $H\!$-type: each group of
$H\!$-type has many different (inequivalent) structures of \ph-type on it.

Section~\ref{cc} contains the formulas for the connection and curvatures,
and gives the explicit forms in each of the examples from
Section~\ref{defex}.  We find a relatively large class of flat spaces, a
clear distinction from the Riemannian case in which there are none.  We
also show that a pseudoeuclidean de Rham factor is characterized in terms
of $j$ when the center is nondegenerate.

Much like the Riemannian case, we would expect that $(N,\langle\, ,\rangle)$
should in some sense be similar to flat pseudoeuclidean space.  This is seen
in the examples of totally geodesic subgroups in Section \ref{geods}.  We
also show the existence of $\dim\z$ independent first integrals, a familiar
result in pseudoeuclidean space, and integrate the geodesic equations, in
certain cases obtaining completely explicit formulas.  Unlike the Riemannian
case, there are flat groups which are isometric to pseudoeuclidean spaces.

Section~\ref{isogp} contains basic information on the isometry group.  In
particular, it can be strictly larger in a significant way than in the
Riemannian (or pseudoriemannian with nondegenerate center) case.

Section~\ref{latpg} begins with a basic treatment of lattices $\G$ in these
groups.  The tori $T_F$ and $T_B$ provide the model fiber and the base for a
submersion of $\G\bs N$.  This submersion may not be pseudoriemannian in the
usual sense, because the tori may be degenerate.  We then begin the study of
periodic geodesics in these compact nilmanifolds, obtaining a complete
calculation of the period spectrum for the flat spaces of Section \ref{cc}.

Finally, section~\ref{cph} details the construction of Lie algebras of
\ph-type from vector spaces with given data (an inner product and a
$j$-operator).  This is a generalization of the construction of Kaplan
\cite{K}, but we have not used the language of Clifford modules here.

We recall some basic facts about 2-step nilpotent Lie groups.  As with all
nilpotent Lie groups, the exponential map $\exp :\n\rightarrow N$ is
surjective.  Indeed, it is a diffeomorphism for simply connected $N$; in
this case we shall denote the inverse by $\log$.  The
Baker-Campbell-Hausdorff formula takes on a particularly simple form in
these groups:
\begin{equation}
\label{bch}
\exp(x)\exp(y)=\exp(x+y+\half[x,y])\, .
\end{equation}
Letting $L_n$ denote left translation by $n\in N$, we have the following
description.
\begin{lemma}
Let\label{dexp} \n\ denote a 2-step nilpotent Lie algebra and $N$ the
corresponding simply connected Lie group. If $x,a\in\n$, then
$$ \exp_{x*}(a_x)=L_{\exp(x)*}\left( a+\half[a,x]\right) $$
where $a_x$ denotes the initial velocity vector of the curve $t\mapsto
x+ta$.\eop
\end{lemma}
\begin{corollary}
In\label{pcc} a pseudoriemannian 2-step nilpotent Lie group, the
exponential map preserves causal character.  Alternatively, 1-parameter
subgroups are curves of constant causal character.
\end{corollary}
\begin{proof}
For the 1-parameter subgroup $c(t) = \exp(t\,a)$, $\dot{c}(t) =
\exp_{ta*}(a) = L_{\exp(ta)*}a$ and left translations are isometries.
\end{proof}
Of course, 1-parameter subgroups need not be geodesics; see Example
\ref{1par}.

We shall also need some basic facts about lattices in $N$. In nilpotent
Lie groups, a lattice is a discrete subgroup \G\ such that the homogeneous
space $M=\G\bs N$ is compact \cite{R}. Lattices do not always exist in 
nilpotent Lie groups \cite{M}.
\begin{theorem}
The simply connected, nilpotent Lie group $N$ admits a lattice if and only 
if there exists a basis of its Lie algebra \n\ for which the structure 
constants are rational.
\end{theorem}
Such a group is said to have a rational structure, or simply to be
rational. 

We recall the result of Marsden from \cite{O}.
\begin{theorem}
A compact, homogeneous pseudoriemannian space is complete.
\end{theorem}
Thus if a rational $N$ is provided with a bi-invariant metric tensor 
$\langle\, ,\rangle$, then $M$ becomes a compact, homogeneous 
pseudoriemannian space which is therefore complete. It follows that 
$(N,\langle\, ,\rangle)$ is itself complete. In general, however, the
metric tensor is not bi-invariant and $N$ need not be complete.

For 2-step nilpotent Lie groups, things work nicely as shown by this result
first published by Guediri \cite{G}.
\begin{theorem}
On\label{gt} a 2-step nilpotent Lie group, all left-invariant pseudoriemannian 
metrics are geodesically complete.
\end{theorem}
\begin{proof}
This follows from the complete integrability of the geodesic equations
(\ref{ges}), or from the integration of them in Theorem \ref{ige}.
\end{proof}
He also provided an explicit example of an incomplete metric on a 3-step
nilpotent Lie group.

\section{Definitions and Examples}
\label{defex}

In the Riemannian (positive-definite)
case, one splits $\n = \z\ds\v = \z\ds\z^{\perp}$ where the
superscript denotes the orthogonal complement with respect to the
inner product $\langle\,,\rangle$. In the general pseudoriemannian
case, however, $\z\ds\z^{\perp}\not=\n$. The problem is that $\z$ might be
a degenerate subspace; {\em i.e.,} it might contain a null subspace $\U$ for
which $\U\subseteq\U^{\perp}$.

Thus we shall have to adopt a more complicated decomposition of \n.
Observe that if $\z$ is degenerate, the null subspace $\U$ is well defined
invariantly. We shall use a decomposition
$$
\n = \z\ds\v = \U\ds\Z\ds\V\ds\E
$$
in which $\z = \U\ds\Z$ and $\v = \V\ds\E$, $\U$ and $\V$ are complementary
null subspaces, and $\U^{\perp}\cap\V^{\perp} = \Z\ds\E$. Although the
choice of $\V$ is {\em not\/} well defined invariantly, once a $\V$ has been 
chosen then $\Z$ and $\E$ {\em are\/} well defined invariantly. Indeed, $\Z$
is the portion of the center $\z$ in $\U^{\perp}\cap\V^{\perp}$ and $\E$
is its orthocomplement in $\U^{\perp}\cap\V^{\perp}$.

We now fix a choice of $\V$ (and therefore $\Z$ and $\E$) to be maintained
throughout this paper. Whenever the effect of this choice is to be
considered, it will be done so explicitly.

Having fixed $\V$, observe that the inner product $\langle\,,\rangle$ provides
a dual pairing between $\U$ and $\V$; {\em i.e.,} isomorphisms $\U^*\cong\V$
and $\U\cong\V^*$. Thus the choice of a basis $\{u_i\}$ in $\U$ determines an
isomorphism $\U\cong\V$ ({\em via\/} the dual basis $\{v_i\}$ in $\V$). 

In addition to the choice of $\V$, we now also fix a basis of $\U$ to be
maintained throughout this paper. Whenever the effect of this choice is to
be considered, it will also be done so explicitly.

We shall also need to use an involution $\io$ that interchanges $\U$ and
$\V$ by this isomorphism and which reduces to the identity on $\Z\oplus\E$
in the Riemannian (positive-definite) case. The choice of such an
involution is not significant, and we have chosen the one which is most
natural for the discussion of $H\!$-type (see {\em infra}). In terms of
chosen orthonormal bases $\{z_\alpha\}$ of $\Z$ and $\{e_a\}$ of $\E$,
$$
\io (u_i)=v_i, \quad \io (v_i)=u_i, \quad \io (z_\alpha )=\ve_\alpha \,
z_\alpha , \quad \io (e_a )= \bar{\ve}_a \, e_a\,,
$$
where, as usual,
$$
\langle u_i,v_i \rangle = 1, \quad \langle z_\alpha ,z_\alpha \rangle =
\ve_\alpha, \quad \langle e_a,e_a \rangle = \bar{\ve}_a\,.
$$
Then $\io (\U )=\V$, $\io (\V )=\U$, $\io (\Z )=\Z$, $\io (\E )=\E$ and
$\io^2 = I$. With respect to this basis $\{u_i,z_\alpha,v_i,e_a\}$ of \n, 
$\io$ is given by the following matrix:
$$
\left[ 
\begin{array}{cccc}
 0 & 0  & \begin{array}{ccc} 1 & \cdots & 0\\[-6pt]
                        \vdots & \ddots & \vdots\\[-3pt]
                             0 & \cdots & 1\end{array} & 0\\
 0 & \begin{array}{ccc}\ve_1 & \cdots & 0\\[-4pt]
                      \vdots & \ddots & \vdots\\[-3pt]
                           0 & \cdots & \ve_r\end{array} & 0 & 0\\
 \begin{array}{ccc} 1 & \cdots & 0\\[-6pt]
               \vdots & \ddots & \vdots\\[-3pt]
                    0 & \cdots & 1\end{array} & 0 & 0 & 0\\
 0 & 0 & 0 & \begin{array}{ccc}\bar{\ve}_1 & \cdots & 0\\[-4pt]
                                    \vdots & \ddots & \vdots\\[-3pt]
                                         0 & \cdots & \bar{\ve}_s\end{array}
\end{array}
\right]
$$
We note that this is also the matrix of $\langle\,,\rangle$ on the same
basis; however, $\io$ is a linear transformation, so it will transform
differently with respect to a change of basis.

It is obvious that $\io$ is self adjoint with respect to the inner product,
\begin{equation}
\langle \io x,y \rangle = \langle x, \io y \rangle , \qquad x,y \in \n\,,
\label{iosa}
\end{equation}
so $\io$ is an isometry of \n. (However, it does not integrate to an
isometry of $N$; see after Example \ref{hqsc}.) Moreover,
\begin{equation}
\langle x, \io x \rangle = 0 \mbox{ if and only if } x = 0, \qquad x \in \n\,.
\label{iosnd}
\end{equation}

Now consider the adjoint with respect to $\langle\,,\rangle$ of the
adjoint representation of the Lie algebra \n\ on itself, to be denoted by
$\add{}{}\!$.  First note that for all $a\in\z$, 
$\add{a}{\!\vcenter{\hbox{$\bul$}}} = 0$.  Thus for all $y\in\n$,
$\add{\bul}{y}$ maps $\V\ds\E$ to $\U\ds\E$.  Moreover, for all $u\in\U$
we have $\add{\bul}{u} = 0$ and for all $e\in\E$ also $\add{\bul}{e} = 0$.
Following \cite{K,E,E'}, we next define the operator $j$. Note the use of
the involution $\io$ to obtain a good analogy to the Riemannian case.
\begin{definition}
The\label{dj} linear mapping
$$ j:\U\ds\Z\rightarrow\End\left(\V\ds\E\right) $$
is given by 
$$ j(a)x = \io\add{x}{\io a}\, .$$
\end{definition}
Equivalently, one has the following characterization:
\begin{equation}\label{cjb}
\langle j(a)x, \io y \rangle = \langle [x, y],\io a \rangle, \quad a\in
\U\ds\Z\, ,\; x,y \in \V\ds\E\, .
\end{equation}
We shall see in Section~\ref{cc} that this map (together with the Lie 
algebra structure of \n) determines the geometry of $N$
just as in the Riemannian case \cite{E'}. 

The operator $j$ is $\io$-skewadjoint with respect to the inner product
$\langle\,,\rangle$. 
\begin{proposition}
For\label{i-sk} every $a\in\U\ds\Z$ and all $x,y\in\V\ds\E$,
$$
\langle j(a)x,\io y \rangle + \langle \io x, j(a)y \rangle = 0\, .
$$
\end{proposition}
\begin{proof}
$$
\langle j(a)x, \io y \rangle = \langle [ x, y],\io a \rangle =
- \langle [ y, x],\io a \rangle = - \langle \io x, j(a)y \rangle \eop
$$
\end{proof}
Although this is not the traditional skewadjoint property, it turns out to
be just what is needed in our generalization of $H\!$-type ({\em
cf.}~Definition \ref{ph} and following). The Riemannian version \cite{K,E'}
is easily obtained as the particular case when we assume $\U=\V=\{0\}$ and
the inner product on $\Z\ds\E$ is positive definite ({\em i.e.,} 
$\ve_\alpha = \bar{\ve}_a=1$); then $\iota =I$ and we recover the
definitions of \cite{K}.  This recovery of the Riemannian case continues
in all that follows.

Recall that the Lie algebra \n\ is said to be nonsingular if and only if
$\ad{x}$ maps \n\ onto $\z$ for every $x\in\n-\z$.
As in \cite{E'}, we immediately obtain
\begin{lemma}
The 2-step nilpotent Lie algebra \n\ is nonsingular if and only if for
every $a\in\U\ds\Z$ the maps $j(a)$ are nonsingular, for every inner
product on \n, every choice of\/ $\V$, every basis of\/ $\U$, and every
choice of $\io$.\eop
\end{lemma}

Next we consider some examples of these Lie groups to which we shall
return in subsequent sections.
\begin{example}
The\label{h3} usual inner products on the 3-dimensional Heisenberg algebra
$\h_3$ may be described as follows. On an orthonormal basis
$\{z,e_1,e_2\}$ the structure equation is
$$ [e_1,e_2]=z $$
with nontrivial inner products
$$ \ve=\langle z,z\rangle ,\quad\bve_1=\langle e_1,e_1\rangle,\quad
\bve_2=\langle e_2,e_2\rangle . $$
We find the nontrivial adjoint maps as
$$
\add{e_1}{z}=\ve\bve_2\, e_2\, ,\quad\add{e_2}{z}=-\ve\bve_1\, e_1
$$
with the rest vanishing. On the basis $\{e_1,e_2\}$,
$$ j(z)=\left[\begin{array}{cc}
            0 & -1\\
            1 & 0
           \end{array}\right] $$
so $j(z)^2=-I_2$. Moreover, a direct computation shows that
$$
j(a)^2 = -\langle a,\io a \rangle I_2 \, , \qquad \mbox{ for all } a\in\Z\,.
$$

This construction extends to the generalized Heisenberg groups $H(p,1)$ of 
dimension $2p+1$ with non-null center of dimension 1 generated by $z$, 
and we find
$$ j(z)=\left[\begin{array}{cc}
0 & -I_p\\
I_p & 0 \end{array}\right]. $$
\end{example}
\begin{example}
In\label{h3n} \cite{CP3} we gave an inner product on the 3-dimensional
Heisenberg algebra $\h_3$ for which the center is degenerate. First we
recall that on an orthonormal basis $\{ e_1,e_2,e_3\}$ with signature
$(+--)$ the structure equations are
\begin{eqnarray*}
[e_3,e_1] &=& \half (e_1-e_2)\, ,\\
{[}e_3,e_2] &=& \half (e_1-e_2)\, ,\\
{[}e_1,e_2] &=& 0\, .
\end{eqnarray*}
The center $\z$ is the span of $e_1-e_2$ and is in fact null.

We take a new basis $\{u=\frac{1}{\sqrt{2}}(e_2-e_1),v=\frac{1}{\sqrt{2}}
(e_2+e_1),e=e_3\}$ and find the structure equation $$[v,e]=u\,.$$ 
Generalizing slightly, we take the nontrivial inner products to be
$$
\langle u,v\rangle = 1\quad\mbox{ and }\quad \langle e,e\rangle = \bve\, .
$$
We find the nontrivial adjoint maps as
$$ \add{v}{v} = \bve e \,,\quad \add{e}{v} = -u\,. $$
On the basis $\{v,e\}$,
$$ j(u) = \left[ \begin{array}{cc}
               0  & -1\\
               1  & 0
              \end{array} \right] $$
so $j(u)^2 = -I_2$. Again, a direct computation shows that
$j(a)^2 = -\langle a,\io a \rangle I_2$ for all $a\in\U$.

Again, this construction extends to $H(p,1)$ with null center generated by $u$
and we find
$$ j(u)=\left[\begin{array}{cc}
                0 & -I_p\\
                I_p & 0 \end{array}\right].$$
\end{example}
According to \cite{CP3}, these are all the Lorentzian inner products on
$\h_3$ up to homothety.
\begin{example}
For\label{hq} the simplest quaternionic Heisenberg algebra of dimension 7, we may
take a basis $\{u_1,u_2,z,v_1,v_2,e_1,e_2\}$ with structure equations
$$ \begin{array}{rclcrcl}
[e_1,e_2]&=&z &\qquad & [v_1,v_2]&=&z\\
{[}e_1,v_1]&=&u_1 && [e_2,v_1]&=&u_2\\
{[}e_1,v_2]&=&u_2 && [e_2,v_2]&=&-u_1
\end{array} $$
and nontrivial inner products
$$ \langle u_i,v_j\rangle=\delta_{ij}\, ,\quad\langle z,z\rangle=\ve\,
,\quad\langle e_a,e_a\rangle=\bve_a\,. $$
As usual, each \ve-symbol is $\pm 1$
independently (this is a combined null and orthonormal basis), so the 
signature is $(++--\,\ve\,\bve_1\,\bve_2)$.

The nontrivial adjoint maps are
$$ \begin{array}{rclcrcl}
\add{v_1}{z}&=&\ve u_2 &\qquad& \add{e_1}{z}&=&\ve\bve_2\,e_2\\[1pt]
\add{v_1}{v_1}&=&-\bve_1\,e_1 && \add{e_1}{v_1}&=&u_1\\[2pt]
\add{v_1}{v_2}&=&-\bve_2\,e_2 && \add{e_1}{v_2}&=&u_2\\[4pt]
\add{v_2}{z}&=&-\ve\,u_1     && \add{e_2}{z}&=&-\ve\bve_1\,e_1\\[1pt]
\add{v_2}{v_1}&=&\bve_2\,e_2  && \add{e_2}{v_1}&=&-u_2\\[2pt]
\add{v_2}{v_2}&=&-\bve_1\,e_1 && \add{e_2}{v_2}&=&u_1
\end{array} $$
For $j$ on the basis $\{v_1,v_2,e_1,e_2\}$ we obtain
$$ j(u_1)=\left[\begin{array}{cccc}
              0  & 0  &  1  & 0\\
              0  & 0  &  0  & -1\\
              -1 & 0  &  0  & 0\\
              0  & 1  &  0  & 0
             \end{array}\right] $$
so $j(u_1)^2=-I_4$,
$$ j(u_2)=\left[\begin{array}{cccc}
              0  &  0 & 0 & 1\\
              0  &  0 & 1 & 0\\
              0  & -1 & 0 & 0\\
              -1 &  0 & 0 & 0
             \end{array}\right] $$
so $j(u_2)^2=-I_4$, 
and
$$ j(z)=\left[\begin{array}{cccc}
            0 & -1 & 0 & 0\\
            1 & 0  & 0 & 0\\
            0 & 0  & 0 & -1\\
            0 & 0  & 1 & 0
           \end{array}\right] $$
so $j(z)^2=-I_4$. Again, a direct computation shows that
$j(a)^2 = -\langle a,\io a \rangle I_4$ for all $a\in\U\ds\Z$.

Once again, this construction may be extended to a nondegenerate center, to 
a degenerate center with a null subspace of dimensions 1 or 3, and to 
quaternionic algebras of all dimensions. 
\end{example}

All of the examples so far are of $H\!$-type; the next one is not.
\begin{example}
For\label{h12} the generalized Heisenberg group $H(1,2)$ of dimension 5 we
take the basis $\{u,z,v,e_1,e_2\}$ with structure equations
\begin{eqnarray*}
[e_1,e_2]&=&z\\
{[}v,e_2]&=&u
\end{eqnarray*}
and nontrivial inner products
$$ \langle u,v\rangle=1\, ,\quad\langle z,z\rangle=\ve\, ,\quad\langle
   e_a,e_a\rangle=\bve_a\,. $$
The signature is $(+-\,\ve\,\bve_1\,\bve_2)$.
We find the nontrivial adjoint maps.
$$ \begin{array}{rclcrcl}
\add{v}{v} &=&\bve_2\,e_2 &\qquad& \add{e_1}{z} &=&\ve\bve_2\,e_2\\[2pt]
\add{e_2}{v} &=&-u       &&       \add{e_2}{z} &=&-\ve\bve_1\,e_1
\end{array} $$
On the basis $\{v,e_1,e_2\}$,
$$ j(u)=\left[\begin{array}{ccc}
            0 & 0 & -1\\
            0 & 0 & 0\\
            1 & 0 & 0
           \end{array}\right] $$
so $j(u)^2\cong 0\ds -I_2$, and
$$ j(z)=\left[\begin{array}{ccc}
            0 & 0 &  0\\
            0 & 0 & -1\\
            0 & 1 & 0
           \end{array}\right] $$
with $j(z)^2=0\ds -I_2$.

This construction, too, extends to the generalized Heisenberg groups $H(1,p)$ 
with $p\ge 2$ which are not of $H\!$-type. The dimension is again $2p+1$,
but now the center has dimension $p$. In each case, the $j$-endomorphisms 
have rank 2 with a similar appearance.
\end{example}

From these examples, it seems natural to extend the usual (positive-
definite) definition of $H\!$-type to indefinite inner products.
\begin{definition}
A\label{ph} $2$-step nilpotent Lie algebra \n\ with indefinite inner product
$\langle\,,
\rangle$ and chosen subspace $\V$ as in Section~\ref{defex}
is said to be of \ph-type (pseudo $H\!$-type) if and only if
$$ j(a)^2 = -\langle a,\io a \rangle I \quad \mbox{ on } \V\ds\E $$
for every choice of $a\in\U\ds\Z$.
\end{definition}
We note that it is easy to verify that this is true for one choice of $\V$
and $\io$ if and only if it is true for any choice of $\V$ and $\io$.  In
the positive-definite case, this reduces to the usual notion of $H\!$-type
\cite{E',K}. In the case of a nondegenerate center, Ciatti \cite{C} has
made an equivalent definition and obtained similar results.  Clearly, if
$\n$ is of \ph-type then $\n$ is nonsingular.

Now it is easy to obtain formulas which are similar to those in
the positive-definite case (see, for example, \cite{E'}, (1.7)).
\begin{proposition}
If\label{phids} $(\n,\langle\,,\rangle,\V)$ is of \ph-type then the following
identities hold:
$$ \begin{array}{lcl}
\langle j(a)x,\io j(a)x \rangle = \langle a,\io a \rangle
\langle x,\io x \rangle & &
\mbox{ for all $a\in\U\ds\Z$,\  $x\in\V\ds\E$}\,;\\
\langle j(a)x,\io j(a)y \rangle = \langle a,\io a \rangle
\langle x,\io y \rangle & &
\mbox{ for all $a\in\U\ds\Z$,\  $x,y\in\V\ds\E$}\,;\\
\langle j(a)x,\io j(b)x \rangle = \langle a,\io b \rangle
\langle x,\io x \rangle & &
\mbox{ for all $a,b\in\U\ds\Z$,\  $x\in\V\ds\E$}\,;\\
j(a) \circ j(b) + j(b) \circ j(a) = -2\langle a,\io b \rangle I& &
\mbox{ for all $a,b\in\U\ds\Z$}\,.
\end{array} $$
\end{proposition}
\begin{proof}
All the identities follow easily taking into account (\ref{iosa}),
Proposition~\ref{i-sk}, and the definition of \ph-type.
 
For the first, we have
$$ \begin{array}{lcl}
\langle j(a)x,\io j(a)x \rangle &=&-\langle \io x,j(a)(j(a)x) \rangle\\
&=&-\langle \io x,-\langle a,\io a \rangle x \rangle\\
&=& \langle a,\io a \rangle \langle x,\io x \rangle .
\end{array} $$

The second and the third identities follow from this one by
polarization.
For the second, we compute as follows:
$$ \langle j(a)(x+y),\io j(a)(x+y) \rangle = \langle a,\io a \rangle
\langle x,\io x \rangle + 2\langle j(a)x,\io j(a)y \rangle +
\langle a,\io a \rangle \langle y,\io y \rangle ; $$
on the other hand,
$$ \begin{array}{lcl}
\langle j(a)(x+y),\io j(a)(x+y) \rangle &=& \langle a,\io a \rangle
\langle x+y,\io (x+y) \rangle \\
&=& \langle a,\io a \rangle \left(\langle x,\io x \rangle + 2\langle x,\io y
\rangle + \langle y,\io y \rangle \right)\,.
\end{array} $$
Similarly for the third:
$$ \langle j(a+b)x,\io j(a+b)x \rangle = \langle a,\io a \rangle
\langle x,\io x \rangle + 2\langle j(a)x,\io j(b)x \rangle +
\langle b,\io b \rangle \langle x,\io x \rangle ; $$
on the other hand,
$$ \begin{array}{lcl}
\langle j(a+b)x,\io j(a+b)x \rangle &=& \langle a+b,\io (a+b) \rangle
\langle x,\io x \rangle \\
&=& \left(\langle a,\io a \rangle + 2\langle a,\io b \rangle
+ \langle b,\io b \rangle \right) \langle x,\io x \rangle .
\end{array} $$
 
In order to prove the fourth identity, we polarize the formula in 
Definition~\ref{ph}:
$$ \begin{array}{lcl}
j(a+b)^2x&=&-\langle a+b,\io (a+b) \rangle x\\
&=&-\left(\langle a,\io a \rangle + 2\langle a,\io b \rangle
+ \langle b,\io b \rangle \right)x\\
&=& j(a)^2x-2\langle a,\io b \rangle x + j(b)^2x\,;
\end{array} $$
on the other hand,
$$ j(a+b)\left(j(a+b)x\right)=j(a)^2x+j(a)\left(j(b)x\right)
+j(b)\left(j(a)x\right)+j(b)^2x\,. \eop $$
\end{proof}

There are no new 2-step nilpotent groups of \ph-type.
\begin{proposition}
The 2-step nilpotent Lie group $N$ is of $H\!$-type if and only if it is
of \ph-type.
\end{proposition}
\begin{proof}
Suppose $(N,\langle\,,\rangle,\V,\io)$ is of \ph-type. Consider the
associated left-invariant Riemannian metric $(\,,)$ given by $(x,y) =
\langle x,\io y\rangle$ and let $j$ and $\tilde{\jmath}$ be the respective
$j$-maps. Then $(\tilde{\jmath}(a)x,y) = ([x,y],a) = \langle[x,y],\io
a\rangle = \langle j(a)x,\io y\rangle$ for all $a\in\z$ and $x,y\in\v$,
noting that $\v$ is the same for both inner products by construction. Thus
$j(a) = \tilde{\jmath}(a)$ for every $a\in\z$, so $\tilde{\jmath}(a)^2 =
j(a)^2 = -\langle a,\io a\rangle I = -(a,a)I$ and $(N,(\,,))$ is of
$H\!$-type.

The converse follows by choosing any $(\,,)$-orthonormal basis, then
choosing an associated $\langle\,,\rangle$ and $\io$, and then reversing
the preceding argument.
\end{proof}
But there are a {\em lot\/} of new geometries now on these groups.
Indeed, given a left-invariant Riemannian metric $(\,,)$ on a Lie group
$N$ of $H\!$-type, the number of nonisometric, associated pseudoriemannian
metrics $\langle\,,\rangle$ grows exponentially with $\dim N$, as is
easily seen by looking at the possible matrices for $\io$.

Thus the notion of $H\!$-type or \ph-type is more fundamentally a group
property than a geometric property. This was already suggested by the
result of Kaplan \cite{K1} that a naturally reductive group of $H\!$-type
has a center of dimension 1 or 3 only.  Our result also shows that
Ciatti's \cite{C} family of admissible algebras with nondegenerate center 
is a subfamily of ours, and thus part of Kaplan's \cite{K}.  See Section
\ref{cph} for a general construction of algebras of \ph-type.

\section{Connection and Curvatures}
\label{cc}

The Levi-Civita connection is given by 
$$
\del{x}{y} = \half\Bigl( [x,y] - \add{x}{y} - \add{y}{x}\Bigr)
$$
for all $x,y\in\n$.
\begin{theorem}
Let\label{conn} $N$ be a 2-step nilpotent Lie group with left-invariant
pseudoriemannian metric tensor $\langle\, ,\rangle$ and Lie algebra $\n =
\U\ds\Z\ds\V\ds\E$ decomposed as in Section~\ref{defex}.
For all $u,u'\in\U$, $z,z'\in\Z$, $v,v'\in\V$, $e,e'\in\E$, and $x\in\V\ds\E$,
we have
\begin{eqnarray}
\del{u}u' = \del{u}z = \del{z}u = \del{z}z' & = & 0\, ,\label{3.1}\\
\del{u}v = \del{v}u = \del{u}e = \del{e}u & = & 0\, ,\\
\del{z}x = \del{x}z &=& -\half \io j(\io z) x\, ,\\
\del{v}v' &=& \half\Bigl( [v,v'] - \io j(\io v) v' - \io j(\io v') 
              v\Bigr)\, ,\label{3.4}\\
\del{v}e &=& \half\Bigl( [v,e] - \io j(\io v) e\Bigr)\, ,\\
\del{e}v &=& \half\Bigl( [e,v] - \io j(\io v) e\Bigr)\, ,\\
\del{e}e' &=& \half [e,e']\,. \eop
\end{eqnarray}
\end{theorem}
\begin{example}
Continuing\label{h3c} from Example~\ref{h3}, the Levi-Civita connection for the
Heisenberg group with non-null center is given by
$$ \begin{array}{rcl}
\del{z}e_1=\del{e_1}z &=& -\half\ve\bve_2\, e_2\\[.5ex]
\del{z}e_2=\del{e_2}z &=& \half\ve\bve_1\, e_1\\
\end{array}\qquad
\del{e_1}e_2=\half z=-\del{e_2}e_1 $$
with the rest vanishing.
\end{example}
\begin{example}
Continuing\label{h3nc} Example~\ref{h3n}, the Levi-Civita connection for the
Heisenberg group with null center is given by
$$ \del{v}v=-\bve\,e \qquad\del{v}e=u $$
with the rest vanishing.
\end{example} 
\begin{example}
Continuing\label{hqc} Example~\ref{hq}, the Levi-Civita connection for our 
quaternionic Heisenberg algebra is given by 
$$ \begin{array}{rclcrcl}
\del{z}v_1=\del{v_1}z&=&-\half\ve\, u_2 &\qquad&
 \del{z}e_1=\del{e_1}z&=&-\half\ve\bve_2\, e_2\\[.5ex]
\del{z}v_2=\del{v_2}z&=&\half\ve u_1 &\qquad&
 \del{z}e_2=\del{e_2}z&=&\half\ve\bve_1\, e_1
\end{array} $$
$$ \begin{array}{rclcrcl}
\del{v_1}v_1&=&\bve_1\, e_1 &\qquad& \del{v_1}e_1&=&-u_1\\
\del{v_1}v_2&=&\half z     &\qquad& \del{v_2}e_1&=&-u_2\\
\del{v_2}v_1&=&-\half z    &\qquad& \del{e_2}v_1&=&u_2\\
\del{v_2}v_2&=&\bve_1\, e_1 &\qquad& \del{e_2}v_2&=&-u_1
\end{array} $$
$$ \del{e_1}e_2= \half z = -\del{e_2}e_1 $$
with the rest vanishing.
\end{example}
\begin{example}
Continuing\label{h12c} Example~\ref{h12}, the Levi-Civita connection for the
group $H(1,2)$ is given by 
$$ \begin{array}{rclcrcl}
\del{z}e_1=\del{e_1}z&=&-\half\ve\bve_2\, e_2 & \qquad &
 \del{z}e_2=\del{e_2}z&=&\half\ve\bve_1\, e_1\\[.5ex]
\del{v}v &=& -\bve_2\, e_2 && \del{v}e_2&=&u
\end{array} $$
$$ \del{e_1}e_2= \half z = -\del{e_2}e_1 $$
with the rest vanishing.
\end{example}

The curvature operator is given by
$$ R(x,y)w = \del{x}\del{y}w - \del{y}\del{x}w - \del{[x,y]}w $$
for all $x,y,w\in\n$.  We denote the projections of vectors in \n\ onto the
subspaces of Section~\ref{defex} by superscripts, $w=w^{\U} + w^{\Z} +
w^{\V} + w^{\E}$, and note that $(\io w)^{\U} = \io w^{\V}$, $(\io w)^{\Z}
= \io w^{\Z}$, $(\io w)^{\V} = \io w^{\U}$, and $(\io w)^{\E} = \io w^{\E}$.
\begin{theorem}
Let\label{curv} $N$ be a 2-step nilpotent Lie group with left-invariant
pseudoriemannian metric tensor $\langle\, ,\rangle$ and Lie algebra $\n =
\U\ds\Z\ds\V\ds\E$ decomposed as in Section~\ref{defex}.  For all
$u,u'\in\U$, $z,z'\in\Z$, $v,v'\in\V$, $e,e'\in\E$, and $x\in\V\ds\E$, we
have
\begin{eqnarray}
R(u,\ )\ &=& R(\ ,\ )u = 0\, ,\\
R(z,z')z'' &=& 0\, ,\\
R(z,z')x   &=& \quar \Bigl(\io j(\io z)[\io j(\io z')x]^{\E} -
               \io j(\io z')[\io j(\io z)x]^{\E}\Bigr)\, ,\\
R(z,x)z'   &=& \quar \io j(\io z)[\io j(\io z')x]^{\E}\, ,\\
R(z,v)v'   &=& \quar\Bigl( [v,\io j(\io z)v']+
               \io j(\io z)[\io j(\io v)v']^{\E}-
               \io j(\io v)[\io j(\io z)v']^{\E}\Bigr. \nonumber\\
           & & {}+\Bigl. \io j(\io z)[\io j(\io v')v]^{\E}\Bigr)\, ,\\
R(z,v)e    &=& \quar\Bigl( [v,\io j(\io z)e] +
               \io j(\io z)[\io j(\io v)e]^{\E}-
               \io j(\io v)[\io j(\io z)e]^{\E} \Bigr)\, ,\\
R(z,e)v    &=& \quar\Bigl( [e,\io j(\io z)v] +
               \io j(\io z)[\io j(\io v)e]^{\E}\Bigr)\, ,\\
R(z,e)e'   &=& \quar [e,\io j(\io z)e']\, ,\\
R(v,v')z   &=& -\quar\Bigl( [\io j(\io z)v,v']+[v,\io j(\io z)v']-
               \io j(\io v)[\io j(\io z)v']^{\E}\Bigr.\nonumber\\
           & & \Bigl.{}+\io j(\io v')[\io j(\io z)v]^{\E}\Bigr)\, ,\\
R(v,v')v'' &=& -\quar\Bigl( [v,\io j(\io v')v''+\io j(\io v'')v']
               -[v',\io j(\io v)v''+\io j(\io v'')v]\Bigr. \nonumber\\
           & &  {}- \io j(\io [v,v'']^{\Z})v'
                +\io j(\io [v',v'']^{\Z})v
                -\io j(\io v)[\io j(\io v')v'']^{\E}  \nonumber\\
           & &  {}+ \io j(\io v')[\io j(\io v)v'']^{\E}
                - \io j(\io v)[\io j(\io v'')v']^{\E} \nonumber\\
           & &  \Bigl.{}+\io j(\io v')[\io j(\io v'')v]^{\E}\Bigr)
                 +\half \io j(\io [v,v']^{\Z})v''\, ,\\
R(v,v')e   &=& -\quar\Bigl( [\io j(\io v)e,v'] +[v,\io j(\io v')e]
                - \io j(\io [v,e]^{\Z})v'\Bigr.\nonumber\\
           & &  \Bigl.{}+\io j(\io [v',e]^{\Z})v 
                -\io j(\io v)[\io j(\io v')e]^{\E} 
                +\io j(\io v')[\io j(\io v)e]^{\E}\Bigr)\nonumber\\
           & &  {}+ \half \io j(\io [v,v']^{\Z})e\, ,\\
R(v,e)z    &=& -\quar\Bigl( [v,\io j(\io z)e] +[\io j(\io z)v,e] -
               \io j(\io v)[\io j(\io z)e]^{\E}\Bigr)\, ,\\
R(v,e)v'   &=& -\quar\Bigl( [v,\io j(\io v')e]+[\io j(\io v)v'+\io j(\io v')
                                                         v,e]\Bigr.\nonumber\\
           & & \Bigl.{}- \io j(\io [v,v']^{\Z})e
               -\io j(\io [v',e]^{\Z})v
               -\io j(\io v)[\io j(\io v')e]^{\E}\Bigr) \nonumber\\
           & & {}+\half \io j(\io [v,e]^{\Z})v'\, ,\\
R(v,e)e'   &=& \quar\Bigl( [e,\io j(\io v)e'] +\io j(\io [v,e']^{\Z})e
                -\io j(\io [e,e']^{\Z})v\Bigr)\nonumber\\
           & & {}+\half \io j(\io [v,e]^{\Z})e'\, ,\\
R(e,e')z   &=& -\quar\Bigl( [e,\io j(\io z)e'] +[\io j(\io z)e,e']\Bigr)\, ,\\
R(e,e')v   &=& -\quar\Bigl( [e,\io j(\io v)e'] +[\io j(\io v)e,e']
               +\io j(\io [v,e]^{\Z})e'\Bigr.\nonumber\\
           & & \Bigl. {}- \io j(\io [v,e']^{\Z})e\Bigr)
               +\half \io j(\io [e,e']^{\Z})v\, ,\\
R(e,e')e'' &=& \quar\Bigl( \io j(\io [e,e'']^{\Z})e'
                          -\io j(\io [e',e'']^{\Z})e\Bigr)
               +\half \io j(\io [e,e']^{\Z})e''\,. \eop 
\end{eqnarray}
\end{theorem}
Together with the formulas for the connection in Theorem~\ref{conn}, these
show that the map $j$ together with the Lie algebra structure of \n\ does
indeed completely determine the geometry of $N$ as in the Riemannian case.
Since $j$ only occurs as $\io j\io = \add{}{}\!$, the appearance of $\io$
is an artifact. It might better be said that the map $\add{}{}\!$ and the
Lie algebra structure of $\n$ completely determine the geometry.

\begin{example}
Continuing\label{h3cc} from Examples \ref{h3} and \ref{h3c}, we find the nontrivial 
curvature operators for the Heisenberg group with non-null center.
$$ \begin{array}{rclcrcl}
R(z,e_1)z&=&-\quar\bve_1\bve_2\,e_1 && R(z,e_2)z&=&-\quar\bve_1\bve_2\,e_2\\
R(z,e_1)e_1&=&\quar\ve\bve_2\, z &\qquad&
 R(e_1,e_2)e_1&=&{\txs\frac{3}{4}}\ve\bve_2\, e_2\\[2pt]
R(z,e_2)e_2&=&\quar\ve\bve_1\, z &\qquad&
 R(e_1,e_2)e_2&=&-{\txs\frac{3}{4}}\ve\bve_1\, e_1
\end{array} $$
\end{example}
\begin{example}
Continuing\label{h3ncc} from Examples \ref{h3n} and \ref{h3nc}, we find that 
all curvatures for the Heisenberg group with null center vanish identically:
this example is flat.

When $p\ge 2$, however, the geometry on $H(p,1)$ is {\em not\/} flat.
\end{example}
Thus there is more than one geometry on the Heisenberg group up to
homothety, in contrast to the Riemannian case of only one.

\begin{example}
Continuing\label{hqcc} from Examples \ref{hq} and \ref{hqc}, we find the 
nontrivial curvature operators for our quaternionic Heisenberg algebra.
$$ \begin{array}{rclcrcl}
R(z,v_1)v_1&=&-\half\ve\bve_1\bve_2\, e_2 &\qquad&
 R(z,e_1)z&=&-\quar\bve_1\bve_2\, e_1\\
R(z,v_1)e_2&=& \half\ve\bve_1\, u_1 &\qquad&
 R(z,e_1)e_1&=&\frac{1}{4}\ve\bve_2\, z\\
R(z,v_2)v_2&=&-\half\ve\bve_1\bve_2\, e_2 &\qquad&
 R(z,e_2)z&=&-\quar\bve_1\bve_2\, e_2\\
R(z,v_2)e_2&=& \half\ve\bve_1\, u_2 &\qquad&
 R(z,e_2)e_2&=&\frac{1}{4}\ve\bve_1\, z\\
&&&&&&\\
R(v_1,v_2)v_1&=&(\bve_1+\frac{3}{4}\ve)u_2 &\qquad&
 R(v_1,e_1)v_2&=&\quar\ve\bve_2\, e_2\\
R(v_1,v_2)v_2&=&-(\bve_1+\frac{3}{4}\ve)u_1 &\qquad&
 R(v_1,e_1)e_2&=& -\quar\ve\, u_2\\[2pt]
R(v_1,v_2)e_1&=& \half\ve\bve_2\, e_2 &\qquad&
 R(v_1,e_2)z&=&-\half\ve\bve_1\, u_1\\
R(v_1,v_2)e_2&=& -\half\ve\bve_1\, e_1 &\qquad&
 R(v_1,e_2)v_1&=& \half\bve_1\, z\\[2pt]
R(v_1,e_2)v_2&=&-\quar\ve\bve_1\, e_1 &\qquad&
 R(v_1,e_2)e_1&=& \quar\ve\, u_2\\
&&&&&&\\
R(v_2,e_1)v_1&=&-\quar\ve\bve_2\, e_2 &\qquad&
 R(v_2,e_1)e_2&=&\quar\ve\, u_1\\[2pt]
R(v_2,e_2)z&=&-\half\ve\bve_1\, u_2 &\qquad&
 R(e_1,e_2)v_1&=&\half\ve\, u_2\\
R(v_2,e_2)v_1&=&\quar\ve\bve_1\, e_1 &\qquad&
 R(e_1,e_2)v_2&=&-\half\ve\, u_1\\[2pt]
R(v_2,e_2)v_2&=&\half\bve_1\, z &\qquad&
 R(e_1,e_2)e_1&=&\frac{3}{4}\ve\bve_2\, e_2\\
R(v_2,e_2)e_1&=&-\quar\ve\, u_1 &\qquad&
 R(e_1,e_2)e_2&=&-\frac{3}{4}\ve\bve_1\, e_1\\
\end{array} $$
\end{example}
\begin{example}
Continuing\label{h12cc} from Examples \ref{h12} and \ref{h12c}, we find the 
nontrivial curvature operators for the group $H(1,2)$.
$$ \begin{array}{rclcrcl}
R(z,v)v&=&-\half\ve\bve_1\bve_2\, e_1 &\qquad& 
 R(z,v)e_1&=&\half\ve\bve_2\, u\\[2pt]
R(z,e_1)z&=&-\quar\bve_1\bve_2\, e_1 &\qquad& 
 R(z,e_1)e_1&=&\quar\ve\bve_2\, z\\[2pt]
R(z,e_2)z&=&-\quar\bve_1\bve_2\, e_2 &\qquad& 
 R(z,e_2)e_2&=&\quar\ve\bve_1\, z\\[2pt]
R(v,e_1)z&=&-\half\ve\bve_2\, u &\qquad&
 R(v,e_1)v&=&\half\bve_2\, z \\
R(e_1,e_2)e_1&=&\frac{3}{4}\ve\bve_2\, e_2 &\qquad&
 R(e_1,e_2)e_2&=&-\frac{3}{4}\ve\bve_1\, e_1
\end{array} $$
\end{example}

Let $x,y\in\n$. Recall that homaloidal planes are those for which the
numerator $\langle R(x,y)y,x\rangle$ of the sectional curvature formula 
vanishes. This notion is useful for degenerate planes tangent to spaces
that are not of constant curvature. 
\begin{theorem}
All\label{zph} central planes are homaloidal: $R(z,z')z'' = R(u,x)y = R(x,y)u 
= 0$ for all $z,z',z''\in\Z$, $u\in\U$, and $x,y\in\n$. Thus the
nondegenerate part of the center is flat:
$$ K(z,z')=0\,. \eop $$
\end{theorem}
This recovers \cite[(2.4),\,item\,c)]{E'} in the Riemannian case. 

In view of this result,
we shall extend the notion of flatness to possibly degenerate submanifolds.
\begin{definition}
A submanifold of a pseudoriemannian manifold is\/ {\em flat}
if and only if every plane tangent to the submanifold is homaloidal.
\end{definition}
\begin{corollary}\label{zf}
The center $Z$ of $N$ is flat.\eop
\end{corollary}
\begin{corollary}\label{ccf}
The only $N$ of constant curvature are flat.
\end{corollary}
\begin{proof}
If $\dim Z > 1$, this follows immediately from the previous Corollary.  If
$\dim Z = 1$ and $\dim N \ge 4$, then from Example \ref{abel} {\em infra\/}
there are abelian subgroups of dimension greater than or equal to 2 which
give rise to flat submanifolds.  Finally, the nonflat Heisenberg group has
nonconstant curvature.
\end{proof}

The degenerate part of the center can have a profound effect on the
geometry of the whole group.
\begin{theorem}
If\/ $[\n,\n] \subseteq \U$ and\label{e0f} $\E=\{0\}$, then $N$ is flat.
\end{theorem}
\begin{proof}
To begin, observe that the first part of the hypothesis implies that $j(z) 
= 0$ for all $z\in\Z$. This eliminates most of the possible contributions 
to a nonzero curvature operator.  With the second part included, one need 
only check (3.17), and this is readily seen to vanish as well.
\end{proof}
Among these spaces, those that also have  $\Z = \{0\}$ (which condition
itself implies $[\n,\n] \subseteq \U$) are fundamental, with the more
general ones obtained by making nondegenerate central extensions.  It is
also easy to see that the product of any flat group with a nondegenerate
abelian factor is still flat.

This is the best possible result in general. Using weaker hypotheses in
place of $\E = \{0\}$, such as $[\V,\V] = \{0\} = [\E,\E]$, it is easy to
construct examples which are not flat.
\begin{corollary}
If\label{exf} $\dim Z \ge \lceil\frac{n}{2}\rceil$, then there exists a flat
metric on $N$.\eop
\end{corollary}
Here $\lceil r\rceil$ denotes the least integer greater than or equal to
$r$.

Before continuing, we pause to collect some facts about the condition
$[\n,\n] \subseteq \U$ and its consequences.
\begin{remarks}{\ns
Since\label{rems} it implies $j(z) = 0$ for all $z\in\Z$, this latter is 
possible with no pseudoeuclidean de Rham factor, unlike the Riemannian
case.

Also, it implies $j(u)$ interchanges $\V$ and $\E$ for all $u\in\U$ if and
only if $[\V,\V] = [\E,\E] = \{0\}$.  Thus an $N$ of \ph-type can have
$[\n,\n] \subseteq \U$ if and only if $\Z = \{0\}$.  Examples are the
Heisenberg group and the groups $H(p,1)$ for $p\ge 2$ with null centers.

Finally we note it implies that, for every $u\in\U$, $j(u)$ maps $\V$ to 
$\V$ if and only if $j(u)$ maps $\E$ to $\E$ if and only if $[\V,\E] = \{0\}$.
}
\end{remarks}

We now continue with some general formulas for sectional curvature.
\begin{theorem}
For\label{kze} any orthonormal $z\in\Z$ and $e\in\E$, 
$$ K(z,e)=\quar\ve_z\bve_e\langle j(\io z)e,j(\io z)e\rangle $$
with $\ve_z=\langle z,z\rangle$ and $\bve_e=\langle e,e\rangle$.\eop
\end{theorem}
This recovers \cite[(2.4),\,item\,b)]{E'} in the Riemannian case.
\begin{theorem}
If\label{kee} $e,e'$ are any orthonormal vectors in $\E$, then 
$$ 
K(e,e')=-{\txs\frac{3}{4}}\bve\bve{\kern.05em}'\langle [e,e'],[e,e']\rangle
$$
with $\bve=\langle e,e\rangle$ and $\bve{\kern.05em}'=\langle e',e'
\rangle$.\eop
\end{theorem}
This recovers \cite[(2.4)\,item\,a)]{E'} in the Riemannian case.

Some other sectional curvature numerators are also relevant.
\begin{proposition}
If\label{oscn} $z\in\Z$, $v\in\V$, and $e\in\E$, then
\begin{eqnarray*}
\langle R(z,v)v,z\rangle &=& \quar\langle j(\io z)v,j(\io z)v\rangle ,\\
\langle R(v,e)e,v\rangle &=& -{\txs\frac{3}{4}}\langle [v,e],[v,e]\rangle
 +\quar\langle j(\io v)e,j(\io v)e\rangle ,\\
\langle R(v,v')v',v\rangle &=& -{\txs\frac{3}{4}}\langle [v,v'],[v,v']\rangle
 +\half\langle j(\io v)v',j(\io v')v\rangle \\
 & & {}+\quar\Bigl( \langle j(\io v')v,j(\io v')v\rangle
     +\langle j(\io v)v',j(\io v)v'\rangle\Bigr)\\
 & & {}-\langle j(\io v)v,j(\io v')v'\rangle . \eop
\end{eqnarray*}
\end{proposition}
\begin{example}
Continuing\label{h3sc} from Examples \ref{h3}, \ref{h3c}, and \ref{h3cc},
we find nontrivial sectional curvatures for the Heisenberg group with
non-null center.
\begin{eqnarray*}
K(z,e_1)=K(z,e_2)&=&\quar\ve\bve_1\bve_2\, ,\\
K(e_1,e_2)&=&-{\txs\frac{3}{4}}\ve\bve_1\bve_2\, .
\end{eqnarray*}
\end{example}
\begin{example}
Continuing\label{hqsc} from Examples \ref{hq}, \ref{hqc}, and \ref{hqcc},
we find sectional curvatures for our quaternionic Heisenberg group.
\begin{eqnarray*}
\langle R(v_1,v_2)v_2,v_1\rangle&=&-(\bve_1+{\txs\frac{3}{4}}\ve)\,,\\
\langle R(v,e)e,v\rangle=\langle R(z,v)v,z\rangle &=& 0\, ,\\
K(z,e_1)=K(z,e_2) &=& \quar\ve\bve_1\bve_2\, ,\\
K(e_1,e_2) &=& -{\txs\frac{3}{4}}\ve\bve_1\bve_2\, .
\end{eqnarray*}
\end{example}
Thus $\io$ cannot integrate to an isometry of $N$ in general, as mentioned
after equation (\ref{iosa}). Isometries must preserve vanishing of
sectional curvature, and an integral of $\io$ would interchange homaloidal
and nonhomaloidal planes in this example.

\begin{example}
Continuing\label{h12sc} from Examples \ref{h12}, \ref{h12c}, and \ref{h12cc},
we find sectional curvatures for the group $H(1,2)$.
\begin{eqnarray*}
\langle R(v,e)e,v\rangle=\langle R(z,v)v,z\rangle &=& 0\, ,\\
K(z,e_1)=K(z,e_2) &=& \quar\ve\bve_1\bve_2\, ,\\
K(e_1,e_2) &=& -{\txs\frac{3}{4}}\ve\bve_1\bve_2\, .
\end{eqnarray*}
\end{example}

The Ricci curvature is a symmetric (0,2)-tensor given by
$$ \Ric(x,y)=\tr\left(\xi\mapsto R(\xi,x)y\right). $$
With respect to the basis $\{u_i,z_{\alpha},v_i,e_a\}$, we have
$$ \Ric(x,y)=\sum_i\langle R(v_i,x)y,u_i\rangle + \sum_{\alpha}\ve_\alpha
\langle R(z_\alpha,x)y,z_\alpha\rangle + \sum_a\bve_a\langle R(e_a,x)y,e_a
\rangle $$
with $\ve_\alpha = \langle z_\alpha,z_\alpha\rangle$ and $\bve_a = \langle
e_a,e_a\rangle$, where we noted that terms of the form $\langle
R(u_i,x)y,v_i \rangle = 0$ for all $i$.
\begin{theorem}
Let\label{ric} $N$ be a 2-step nilpotent Lie group with left-invariant
pseudoriemannian metric tensor $\langle\, ,\rangle$ and Lie algebra $\n =
\U\ds\Z\ds\V\ds\E$ decomposed as in Section~\ref{defex}.
For all $u\in\U$, $z,z'\in\Z$, $v,v'\in\V$, $e\in\E$, and $x\in\V\ds\E$,
we have
\begin{eqnarray}
\Ric(u,\ ) &=& 0\,,\\
\Ric(z,z') &=& \quar\sum_a\bve_a\langle j(\io z)e_a,j(\io z')e_a\rangle ,\\
\Ric(v,z) &=& \quar\sum_a\bve_a\langle j(\io v)e_a,j(\io z)e_a\rangle ,\\
\Ric(e,z) &=& 0\,,\\[.5ex]
\Ric(v,v') &=& -\half\sum_\alpha\ve_\alpha\langle j(z_\alpha)v,j(z_\alpha)v'
  \rangle + \quar\sum_a\bve_a\langle j(\io v)e_a,j(\io v')e_a\rangle ,
   \nonumber\\[-2.5ex]
          &&\\
\Ric(x,e) &=& -\half\sum_\alpha\ve_\alpha\langle j(z_\alpha)x,j(z_\alpha)e
  \rangle .
\end{eqnarray}
\end{theorem}
\begin{proof}
This is mostly just straight-forward computation, but there are some
details to be noted carefully. We give one case to illustrate this, in
which we use the selfadjointness of $\io$ and the $\io$-skewsymmetry of
$j$ extensively. From the definition,
$$ \Ric(x,e) = \sum_\alpha\ve_\alpha\langle R(x,z_\alpha)z_\alpha,e\rangle
+\sum_a\bve_a\langle R(x,e_a)e_a,e\rangle . $$
From Theorem \ref{curv} and those properties of $\io$ and $j$, 
\begin{eqnarray*}
\langle R(x,z_\alpha)z_\alpha,e\rangle &=& -\quar\langle\io j(\io
z_\alpha)[\io j(\io z_\alpha)x]^{\E},e\rangle\\
&=& \quar\langle j(z_\alpha)x,j(z_\alpha)e\rangle .
\end{eqnarray*}
The basic trick in all these cases is to expand in an orthonormal basis and 
write
$$ j\left(\io[e_a,x]^{\Z}\right) = \sum_\alpha\ve_\alpha\langle\io j(\io
z_\alpha)e_a,x\rangle j(\io z_\alpha)\,. $$
Then we get
\begin{eqnarray*}
\langle R(e_a,x)e_a,e\rangle &=& {\txs\frac{3}{4}}\sum_\alpha\langle\io
  j(\io z_\alpha)e_a,x\rangle\langle\io j(\io z_\alpha)e_a,e\rangle\\
&=& {\txs\frac{3}{4}}\sum_\alpha\ve_\alpha\langle j(\io
  z_\alpha)x,e_a\rangle \langle j(\io z_\alpha)e,e_a\rangle\\
&=& {\txs\frac{3}{4}}\sum_\alpha\ve_\alpha\langle j(z_\alpha)x,e_a\rangle 
  \langle j( z_\alpha)e,e_a\rangle ,
\end{eqnarray*}
and summing on $a$ yields
$$ \sum_a\bve_a\langle R(e_a,x)e_a,e\rangle = {\txs\frac{3}{4}}\sum_\alpha
  \ve_\alpha\langle j(z_\alpha)x,j( z_\alpha)e\rangle . $$
This results in the formula for $\Ric(x,e)$, and the others are similar.
\end{proof}
\begin{corollary}
If\label{ricf} $j(z) = 0$ for all $z\in\Z$ and $j(u)$ interchanges $\V$
and $\E$ for all $u\in\U$, then $N$ is Ricci flat.\eop
\end{corollary}
\begin{theorem}
Let\label{scal} $N$ be a 2-step nilpotent Lie group with left-invariant
pseudoriemannian metric tensor $\langle\, ,\rangle$ and Lie algebra $\n =
\U\ds\Z\ds\V\ds\E$ decomposed as in Section~\ref{defex}.
The scalar curvature is given by
$$ S = -\quar\sum_{\alpha,a}\ve_\alpha\bve_a\langle j(z_\alpha)e_a,
    j(z_\alpha)e_a\rangle . \eop $$
\end{theorem}
\begin{corollary}
If\label{scalf} $j(z) = 0$ for all $z\in\Z$, then $N$ is scalar flat. In
particular, this occurs when $[\n,\n]\subseteq\U$.
\end{corollary}
\begin{proof}
The second part follows from the observation at the beginning of the proof
of Theorem \ref{e0f}.
\end{proof}
\begin{example}
Continuing\label{h3rics} from Examples \ref{h3}, \ref{h3c}, \ref{h3cc}, and 
\ref{h3sc}, we find nontrivial Ricci curvatures for the Heisenberg group 
with non-null center.
\begin{eqnarray*}
\Ric(z,z) &=& \half\bve_1\bve_2\\
\Ric(e_1,e_1) &=& -\half\ve\bve_2\\
\Ric(e_2,e_2) &=& -\half\ve\bve_1
\end{eqnarray*}
The scalar curvature is $S = -\half\ve\bve_1\bve_2$.

For the group $H(p,1)$ with non-null center, we take a basis of the Lie
algebra $\{z,e_1,\ldots,e_p,e_{p+1},\ldots,e_{2p}\}$ with structure
equation $[e_i,e_{p+i}] = z$ for $1\le i\le p$. Then we find
\begin{eqnarray*}
\Ric(z,z) &=& \half\sum_{i=1}^p\bve_i\bve_{p+i}\,,\\
\Ric(e_i,e_i) &=& -\half\ve\bve_{p+i}\,,\\
\Ric(e_{p+i},e_{p+i}) &=& -\half\ve\bve_i\,.
\end{eqnarray*}
The scalar curvature is $S = -\half\ve\sum_i\bve_i\bve_{p+i}$.
\end{example}
\begin{example}
Continuing\label{h3nrics} from Examples \ref{h3n}, \ref{h3nc}, and 
\ref{h3ncc}, since the Hei\-senberg group with null center is flat, it is
also Ricci flat and scalar flat.

For the group $H(p,1)$ with $p\ge2$ and null center, we take a basis of 
the Lie algebra $\{u,v,e_2,\ldots,e_p,e_{p+1},\ldots,e_{2p}\}$ with structure
equations $[v,e_{p+1}] = u$ and $[e_i,e_{p+i}] = u$ for $2\le i\le p$. 
Then we find
$$ \Ric(v,v) = \half\sum_{i=2}^p\bve_i\bve_{p+i}\,, $$
but the scalar curvature vanishes: $S=0$.  Note that for $p$ odd, with an
appropriate choice of signature one obtains examples of metrics that are
not flat, but are Ricci and scalar flat.
\end{example}
\begin{example}
Continuing\label{hqrics} from Examples \ref{hq}, \ref{hqc}, \ref{hqcc}, and 
\ref{hqsc}, we find nontrivial Ricci curvatures for our quaternionic 
Heisenberg algebra.
\begin{eqnarray*}
\Ric(z,z) &=& \half\bve_1\bve_2\\
\Ric(e_1,e_1) &=& -\half\ve\bve_2\\
\Ric(e_2,e_2) &=& -\half\ve\bve_1
\end{eqnarray*}
The scalar curvature is $S = -\half\ve\bve_1\bve_2$.
\end{example}
\begin{example}
Continuing\label{h12rics} from Examples \ref{h12}, \ref{h12c}, \ref{h12cc},
and \ref{h12sc}, we find nontrivial Ricci curvatures for the group
$H(1,2)$.
\begin{eqnarray*}
\Ric(z,z) &=& \half\bve_1\bve_2\\
\Ric(e_1,e_1) &=& -\half\ve\bve_2\\
\Ric(e_2,e_2) &=& -\half\ve\bve_1
\end{eqnarray*}
The scalar curvature is $S = -\half\ve\bve_1\bve_2$.

For the group $H(1,p)$ with $p\ge2$, we take a basis $\{u_1,\ldots,u_q,
z_1,\ldots,z_{p-q},$\newline$ v_1,\ldots,v_q, e_1,\ldots,e_{p-q}, e\}$ with
structure equations $[v_i,e] = u_i$ for $1\le i\le q$ and $[e_j,e] = z_j$
for $1\le j\le p-q$. Then we find
\begin{eqnarray*}
\Ric(z_j,z_j) &=& \half\bve_j\bve\,,\\
\Ric(e_j,e_j) &=& -\half\ve_j\bve\,,\\
\Ric(e,e) &=& -\half\sum_{j=1}^{p-q}\ve_j\bve_j\,.
\end{eqnarray*}
The scalar curvature is $S = -\half\bve\sum_j\ve_j\bve_j$.
\end{example}
We note the amusing low-dimensional coincidence of three of these examples.

When $N$ has a pseudoeuclidean de Rham factor, it may be characterized in
terms of $j$; {\em cf.} \cite[(2.7)]{E'}.  To use the pseudoriemannian de
Rham theorem of Wu \cite{W}, we must assume that $\z$ is nondegenerate.
\begin{proposition}
Let $N$ be a simply connected, 2-step nilpotent Lie group with
left-invariant metric tensor $\langle\,,\rangle$ and nondegenerate center.
Assume that $N$ can be written as a pseudoriemannian product $E\times N'$
with $E$ the pseudoeuclidean de Rham factor and $N'$ the product of all
other de Rham factors, and write $\n = \e\ds\n'$ as the corresponding
orthogonal direct sum. Then $\e\subseteq\z$ and $\e = \{z\in\z\mid
j(z)=j(\io z)=0\}=\ker j\cap\ker j\io$.
\end{proposition}
\begin{proof}
For convenience, we shall also denote the corresponding subbundles of $TN$
by $\e$ and $\n'$. Then we note that $\e$ is parallel, integrable, and
$\io$-invariant. Also, $\n=\z\ds\v$ is an orthogonal direct sum and $\z$
and $\v$ are $\io$-invariant, since we have assumed $\z$ is nondegenerate.

First we show that $\e\cap\v=\{0\}$. Suppose not, with $0\ne
x\in\e\cap\v$. Then for all $z\in\z$, we have $j(z)x = -2\io\del{\io z}x$.
Since $\e$ and $\v$ are $\io$-invariant, $\del{\io z}x\in\v$, and $\e$ is
parallel, we deduce that $j(z)x\in\e\cap\v$ for all $z\in\z$. To see that
$j(z)x\ne 0$ for some $z\in\z$, choose $y\in\v$ such that $[x,y]\ne 0$.
Then for $z=[x,y]$, we find $\langle\io j(z)x,y\rangle = \langle\add{x}\io
z ,y\rangle = \langle[x,y],\io z\rangle = \langle z,\io z\rangle\ne 0$.
Thus $\io j(z)x\ne 0$ whence 
\begin{equation}
\label{bbb}
0\ne x\in\v\mbox{ implies }j(z)x\ne 0\mbox{ for some }z\in\z\,. 
\end{equation}
Now, for any such $z$, we have
$[x,j(z)x]\in\e\cap\z$ because $x\in\e\cap\v$, $j(z)x\in\e\cap\v$, and
$\e$ is integrable. Moreover, for such a $z$, we would have $\langle[x,
j(z)x],\io z\rangle = \langle\io j(z)x,j(z)x\rangle\ne 0$ and $[x,j(z)x]$
would be another nonzero element in $\e\cap\v$. But this would mean that
a sectional curvature numerator for $\lsp x,j(z)x\rsp$, namely
$\langle[x,j(z)x],[x,j(z)x]\rangle$, is nonvanishing. This contradicts $E$
being flat, so we conclude that $\e\cap\v=\{0\}$ as desired.

Next we show that $\e\subseteq\z$. Let $0\ne a\in\e$ and write $a=z+x$
with $z\in\z$ and $x\in\v$. To show that $x=0$ it suffices by (\ref{bbb})
to show that $j(z')x=0$ for all $z'\in\z$. Let $z'\in\z$ be fixed but
arbitrary. Then the parallel and $\io$-invariant distribution $\e$
contains $\del{\io z'}a = \del{\io z'}x = -\half\io j(z')x$. But from the
first part of this proof it follows that $j(z')x=0$ since
$j(z')x\in\e\cap\v$. Therefore $\e\subseteq\z$.

We now claim that for all $z\in\e$, we have $j(z)=j(\io z)=0\in\End(\v)$.
Indeed, for any $z\in\e$ and $x\in\v$, since $\e$ is parallel and
$\io$-invariant, $\io\del{x}z=-\half j(\io z)x$. But $j(\io
z)x\in\e\cap\v=\{0\}$; similarly, $j(z)x=0$.

Finally, we show that if $j(z)=j(\io z)=0$ then $z\in\e$. For any such $z$
set $\e'=\lsp\e,z\rsp$, and observe that it suffices to prove that $\e'$
is parallel and flat. Consider $a'\in\e'$ and $x'\in\n$ and write $a'=a+tz$
for some $a\in\e$ and $t\in\R$, and $x'=z'+x$ for some $z'\in\z$ and
$x\in\v$. Then $$\del{x'}a' = \del{x'}a+t\del{x'}z =
\del{x'}a-{\txs\frac{t}{2}}\io j(\io z)x = \del{x'}a\in\e\subseteq\e'$$ so
$\e'$ is parallel. Now, from $\e'\subseteq\z$ it follows that
$\del{e}e'=0$ for all $e,e'\in\e'$, whence $\e'$ is flat.
\end{proof}

\section{Totally Geodesic Subgroups and Geodesics}
\label{geods}

We begin by noting that O'Neill \cite[Ex.\,9, p.\,125]{O} has extended the
definition of totally geodesic to degenerate submanifolds of
pseudoriemannian manifolds. We shall use this extended version.

Recall from \cite{O} that the extrinsic and intrinsic curvatures of
totally geo\-des\-ic submanifolds coincide. Thus there is an unambiguous
notion of flatness for them.  Note that a connected subgroup $N'$ of $N$
is a totally geodesic submanifold if and only if it is totally geodesic at
the identity element of $N$, because left translations by elements of $N'$
are isometries of $N$ that leave $N'$ invariant. A connected, totally
geodesic submanifold need not be a connected, totally geodesic subgroup,
but $(N,\langle\, ,\rangle)$ has many totally geodesic subgroups. Many of
them are flat, illustrating the similarity to pseudoeuclidean spaces; {\em
cf.} \cite[(2.11)]{E'}.
\begin{example}
For\label{1par} any $x\in\n$ the 1-parameter subgroup $\exp (tx)$ is a
geo\-des\-ic if and only if $\add{x}{x}=0$. We find this if and only if
$x\in\z$ or $x\in\U\ds\E$. This is essentially the same as the Riemannian
case, but with some additional geodesic 1-parameter subgroups coming from
$\U$.
\end{example}
\begin{example}
Abelian\label{abel} subspaces of $\V\ds\E$ are Lie subalgebras of \n, and
give rise to complete, flat, totally geodesic abelian subgroups of $N$,
just as in the Riemannian case \cite{E'}. The construction given in
\cite[(2.11),\,Ex.\,2]{E'} is valid in general, and shows that if $\dim\V
\ds\E \geq 1 + k + k\dim\z$, then every nonzero element of $\V\ds\E$ lies
in an abelian subspace of dimension $k+1$.
\end{example}
\begin{example}
The\label{ftgZ} center $Z$ of $N$ is a complete, flat, totally geodesic
submanifold.  Moreover, it determines a foliation of $N$ by its left
translates, so each leaf is flat and totally geodesic, as in the
Riemannian case \cite{E'}. In the pseudoriemannian case, this foliation in
turn is the orthogonal direct sum of two foliations determined by $\U$ and
$\Z$, and the leaves of the $\U$-foliation are also null.  All these
leaves are complete.
\end{example}

Let $\g$ be a geodesic in a group $N$ of \ph-type. Without loss of
generality, we may assume that $\g(0) = 1\in N$. Write $\dot{\g}(0) = a_0
+ x_0 \in(\U\ds\Z)\ds(\V\ds\E)$ and suppose that neither $a_0$ nor $x_0$
is zero. Consider the span $\lsp a_0,x_0,j(a_0)x_0\rsp = \n'$ in the Lie
algebra \n\ of $N$. The following lemma shows that $\n'$ is isomorphic to
the 3-dimensional Heisenberg algebra $\h_3$.
\begin{lemma}
Let $\n = \U\ds\Z\ds\V\ds\E$ decomposed as in Section~\ref{defex} be an
algebra of \ph-type.
If $x\in\V\ds\E$ and $a\in\U\ds\Z$, then $[x,j(a)x] = \langle x,\io
x\rangle a$.
\end{lemma}
\begin{proof}
For any $a'\in\U\ds\Z$, using (\ref{cjb}) and the third identity from
Proposition~\ref{phids}, we have
\begin{eqnarray*}
\langle[x,j(a)x],\io a'\rangle &=& \langle j(a)x,\io j(a')x\rangle\\
&=& \langle x,\io x\rangle\langle a,\io a'\rangle .
\end{eqnarray*}
Thus $\langle a,\io a'\rangle = 0$ implies there is no component in the
$a'$ direction.
\end{proof}
In the Riemannian case, each such subgroup $N'$ with Lie algebra $\n'$ is
totally geodesic \cite[p.\,624f\,]{E'}. But in general, this is not true.
\begin{example}
In our quaternionic Heisenberg group from Example~\ref{hq}, consider a
geodesic $\g$ with $\dot{\g}(0) = a_0 + x_0 = (u_2+z)+(v_1+e_1)$.  Then
one computes easily $j(a_0)x_0 = 2v_2$, so from Example~\ref{hqc} we find
$\del{x_0}a_0 = -\half\ve(u_2+\bve_2 e_2)\notin\n' = \lsp u_2+z, v_1+e_1,
v_2\rsp$.
\end{example}
As usual, part of the problem is the presence of the null subspaces $\U$ 
and $\V$. But even if we assume that $\z$ is nondegenerate, these subgroups
still need not be totally geodesic.
\begin{example}
For the group $H(p,1)$ with $p\ge 2$, take $a=z$ and $x=e_i+e_j$ with $i\ne
j$. Then $j(a)x = j(z)(e_i+e_j) = e_{p+i}+e_{p+j}$ and 
$$ \del{a}x = -\half\ve(\bve_{p+i}e_{p+i} + 
\bve_{p+j}e_{p+j}) \notin \lsp a,x,j(a)x\rsp\,. $$
One may also check that using, for example, $\io j(\io a)x=\add{x}a$ will
not help: in the quaternionic Heisenberg algebra with $a=z_1$ and $x =
e_1+e_4$, $\lsp a,x,\add{x}a\rsp$ is totally geodesic but is not a
subalgebra.
\end{example}

For the geodesic equation, let us consider (as suffices) a geodesic $\g$
with $\g(0) = 1\in N$ and $\dot{\g}(0) = a_0+x_0\in \z\ds\v$. Further
decompose $a_0 = u_0 + z_0 \in \U\ds\Z$ and $x_0 = v_0 + e_0 \in \V\ds\E$.
In exponential coordinates, write $\g(t) = \exp \left( a(t)+x(t)\right) = 
\exp\left( u(t) + z(t) + v(t) + e(t)\right)$ with $\dot{a}(0) = a_0$,
$\dot{x}(0) = x_0$, {\em etc.} For the tangent vector $\dot{\g}$, we obtain
\begin{eqnarray*}
\dg &=& \exp_{(a+x)*}(\dot{a}+\dx)\\
   &=& L_{\g(t)*}\left(\dot{a}+\dx + \half[\dot{a}+\dx,a+x]\right)\\
   &=& L_{\g(t)*}\left(\dot{a}+\dx + \half[\dx,x]\right),
\end{eqnarray*}
using Lemma \ref{dexp},
regarded as vector fields along $\g$. Then the geodesic equation is
equivalent to
$$ \frac{d}{dt}\left( \dot{a}+\dx + \half[\dx,x]\right)
   - \add{\dot{a} + \dx +\frac{1}{2}[\dx,x]}{\left(\dot{a}
   + \dx +\half[\dx,x]\right)} = 0\,. $$
Simplifying slightly, we find
\begin{equation}
\underbrace{\frac{d}{dt}\left(\dot{a} + \half[\dx,x]
\right)_{\strut}}_{\in\,\U\ds\Z}{\;} +
\underbrace{\ddx_{\strut}}_{\in\,\V\ds\E} -
{\;\,}\underbrace{\add{\dx}{\left(\dot{z} + \dot{v} +
\half[\dx,x]\right)}_{\strut}}_{\in\,\U\ds\E} = 0\,.\label{ges}
\end{equation}
Using superscripts to denote components as before, we obtain for the
$\Z$-component
$$ \frac{d}{dt}\left(\dot{z} + \half[\dx,x]^{\Z}\right) = 0\,. $$
Using the initial condition, we get
$$ \dot{z} + \half[\dx,x]^{\Z} = z_0\,. $$
Next we note that $\ddot{v} = 0$ whence $\dot{v}(t) = v_0$ is a
constant. We use these to simplify the other component equations.
\begin{eqnarray}
\frac{d}{dt}\left(\dot{u} + \half[\dx,x]^{\U}\right) -
   \left(\add{\dx}{(z_0+v_0)}\right)^{\U} &=& 0\label{ges1}\\
\dot{z} + \half[\dx,x]^{\Z} &=& z_0\label{ges2}\\
\ddot{v} &=& 0\label{ges3}\\
\ddot{e} - \left(\add{\dx}{(z_0+v_0)}\right)^{\E} &=& 0\label{ges4}
\end{eqnarray}

In analogy with Eberlein \cite{E,E'} we define two operators.
\begin{definition}
For\label{defJ} fixed $z_0\in\Z$ and $v_0\in\V$ as above, define
\begin{eqnarray*}
\rsfs{J}:\V\ds\E \longrightarrow \E &:& y \longmapsto 
\left(\add{y}{(z_0+v_0)}\right)^{\E},\\
\sJ:\V\ds\E \longrightarrow \U &:& y \longmapsto 
\left(\add{y}{(z_0+v_0)}\right)^{\U}.
\end{eqnarray*}
\end{definition}
We shall denote the restriction of $\rsfs{J}$ to $\E$ by $J$, and this will
play the same role as $J$ in Eberlein \cite{E,E'}.

Now we rewrite the geodesic equations in terms of $J$ and $\sJ$,
using the linearity of $\add{}{}$ to rearrange some terms.
\begin{eqnarray}
\frac{d}{dt}(\dot{u} + \half[\dx,x]^{\U}) &=& \sJ\dx \label{ge1}\\
\dot{z} + \half[\dx,x]^{\Z} &=& z_0\label{ge2}\\
\ddot{v} &=& 0\label{ge3}\\
\ddot{e} - J\dot{e} &=& \rsfs{J}v_0\label{ge4}
\end{eqnarray}
While the $\V$-component of a geodesic is simple, its mere presence
affects all of the other components.

We also readily see that the system is completely integrable.  Thus, as
noted in the Introduction (proof of Theorem \ref{gt}), all left-invariant
pseudoriemannian metrics on these groups are complete.  Also, {\em
regardless of signature,} we may obtain the existence of $\dim\z$ first
integrals as in \cite{E'}.

Keep $J\in\End(\E)$ and write $\E = \E_1\ds\E_2$ with $\E_1 = \ker J$, an
orthogonal direct sum. Decompose $\rsfs{J}v_0 = y_1 + y_2 \in
\E_1\ds\E_2$, respectively.  Note that $J$ is invertible on $\E_2$; we
denote this restriction by $J$ also.  Now we follow Eberlein
\cite[(3.2--3.4)]{E'} for these next three results.
\begin{proposition}
Let\label{e3.2} $N$ be simply connected, $\g$ a geodesic with $\g(0)=1\in N$,
and $\dot{\g}(0) = a_0 + x_0 \in \n$. Then
$$\dot{\g} = L_{\gamma(t)*}\left( a_0 + e^{tJ}(e_0+J^{-1}y_2) + v_0 -
J^{-1}y_2 + t\,y_1 + \I\sJ\dx\right) $$
where $\I\sJ\dx = \int_0^t\sJ\dx(s)\,ds$.\eop
\end{proposition}
\begin{corollary}
Let\label{e3.3} \F\ denote the geodesic flow in $TN$, $n\in N$, and $a_0+x_0
\in \n$. Then
$$\F L_{n*}(a_0+x_0) = L_{\gamma(t)*}\left(a_0 + e^{tJ}(e_0+J^{-1}y_2) +
v_0 - J^{-1}y_2 + t\,y_1 + \I\sJ\dx\right) $$
where $\g$ is the unique geodesic with $\g(0)=n$ and $\dot{\g}(0) = 
L_{n*}(a_0 + x_0)$.\eop
\end{corollary}
The proofs are direct computations which we omit; see \cite{E,E'}.

In the next result, we use $PTM$ to denote the projectivized tangent
bundle of a manifold $M$.
\begin{corollary}
Let $\{z_1,\ldots ,z_r\}$ be any basis of $\z$ and
let $\{z^*_1,\ldots ,z^*_r\}$ be the dual
basis of $\z^*$. Define $f_0 : TN\rightarrow \z : L_{n*}x\mapsto \pi x$
where $n\in N$, $x\in\n$, and $\pi : \n\rightarrow\z$ is the projection.
For $1\le\a\le r$, define $f_\a : TN\rightarrow\R$ by
$$
f_\a= z^*_\a \circ f_0\,, \quad 1\le \a\le r\, .
$$
Then for $0\le \a\le r$, $t\in\R$, and $n\in N$, we have these:
\begin{enumerate}
\item $f_\a\circ\F = f_\a$.
\item $f_\a\circ L_{n*} = f_\a$.
\item If\/ $\G\le N$ is any discrete subgroup acting on $N$ by left
translations, then $f_\a$ induces $F_\a : T(\G\bs N)\rightarrow\R$ with $F_\a
\circ\F = F_\a$ for \F\ the induced geodesic flow on $T(\G\bs N)$. In
particular, \F\ has no dense orbit in any of the three parts (timelike,
spacelike, or null) of $PT(\G\bs N)$.
\item The functions $f_\a$ and $F_\a$ for $1\le \a \le r$ are linearly
independent over \R.
\end{enumerate}
\end{corollary}
\begin{proof}
The first two parts are immediate consequences of Corollary \ref{e3.3}.
Letting $p_N:N\surj\G\bs N$ denote the natural projection, the second part
implies that the functions $F_0:T(\G\bs N)\to\z$ and $F_\a:T(\G\bs N)\to\R$
for $\a\ge1$
are well defined and satisfy $F_\a p_{N*}=f_\a$ for all $\a$. From the
first part, the $F_\a$ commute with \F. The functions $F_\a$ are
clearly continuous and nonconstant in $PT(\G\bs N)$, but constant along
\F-orbits. Thus there is no dense orbit in any part of $PT(\G\bs N)$.
Finally, the last part follows from the definition of $F_\a$ and $f_\a$.
\end{proof}

We continue with a geodesic through the identity element, so $\g(0)=1$ and
$\dot{\g}(0) = a_0+x_0 = u_0+z_0+v_0+e_0$, with $J\in\End(\E)$ and $\E =
\E_1\ds\E_2$ with $\E_1 = \ker J$, and with $\rsfs{J}v_0 = y_1 + y_2 \in 
\E_1\ds\E_2$ as before.  Decompose $e_0 = e_1 + e_2 \in \E_1\ds\E_2$,
respectively.  (These $e_i$ should not be confused with the basis elements
appearing in other sections.)  Observe that $J$ is skewadjoint with
respect to $\langle\,,\rangle$.  Recall that $J$ is invertible on $\E_2$,
and that we let $J$ denote the restriction there as well.  For
convenience, set
\begin{eqnarray*}
x_1 &=& e_1 + v_0 - J^{-1}y_2\,,\\
x_2 &=& e_2 + J^{-1}y_2\,.
\end{eqnarray*}
\begin{theorem}
Using\label{ige} these notations, the geodesic equations may be integrated
as:
\begin{eqnarray}
x(t) &=& t\,x_1 + \left( e^{tJ} - I\right) J^{-1}x_2 + \half t^2y_1\,,
         \label{ige1}\\
z(t) &=& t\,z_0 + \I[\dx,x]^{\Z},\label{ige2}\\
u(t) &=& t\,u_0 + \I[\dx,x]^{\U} + \II\sJ\dx\,,\label{ige3}
\end{eqnarray}
where
\begin{eqnarray*}
 \I[\dx,x] &=& -\half\int_0^t\left[\dx(s),x(s)\right] ds\\
&=&\half t\left[ x_1 + \half ty_1, \left(e^{tJ} + I\right) J^{-1}x_2 \right]
- t\left[ y_1, e^{tJ}x_2\right] + {\txs\frac{1}{12}}t^3[x_1, y_1]\\
&&{}+ \half\left[\left( e^{tJ} - I\right) J^{-1}x_2, J^{-1}x_2\right]
- \left[ x_1, \left( e^{tJ} - I\right) J^{-2}x_2\right]\\
&&{}+ \left[ y_1, \left( e^{tJ} - I\right) J^{-3}x_2\right]
+ \half\int_0^t \left[ e^{sJ} J^{-1}x_2, e^{sJ} x_2\right] ds
\end{eqnarray*}
and
$$ \II\sJ\dx = \int_0^t\int_0^s\sJ\dx(\sigma)\,d\sigma\,ds\,. $$
\end{theorem}
\begin{proof}
The formulas follow from straightforward integrations of the geodesic
equations (\ref{ge1})--(\ref{ge4}).  We used the general fact about
exponentials of matrices that $J$ commutes with $e^{tJ}$ for all $t\in\R$.
Using this, it is routine to verify that $x(t)$, $z(t)$, and $u(t)$
satisfy the geodesic equations and initial conditions.
\end{proof}
\begin{corollary}
When\label{igend} $N$ has a nondegenerate center, the formulas simplify
somewhat.  Now equation (\ref{ge4}) is homogeneous and we obtain 
\begin{eqnarray}
e(t) &=& t\,e_1 + \left( e^{tJ} - I\right) J^{-1}e_2\,,\label{igend1}\\
z(t) &=& t\,z_1(t) + z_2(t) + z_3(t)\,,\label{igend2}
\end{eqnarray}
where
\begin{eqnarray*}
z_1(t) &=& z_0+\half\left[e_1,\left(e^{tJ}+I\right)J^{-1}e_2\right],\\
z_2(t) &=& \left[e_1,\left(I-e^{tJ}\right) J^{-2}e_2\right] + 
           \half\left[e^{tJ}J^{-1}e_2, J^{-1}e_2\right],\\
z_3(t) &=& \half\int_0^t \left[e^{sJ}J^{-1}e_2, e^{sJ}e_2\right]
      \,ds\,. 
\end{eqnarray*}
Note that $z_3$ may contribute to $z_1$ and $z_2$.\eop
\end{corollary}
The flat spaces we found in Theorem \ref{e0f} allow more simplification.
\begin{corollary}
When\label{iged} $[\n,\n] \subseteq \U$ and $\E = \{0\}$, then $x = v$ and
we obtain
\begin{eqnarray}
v(t) &=& t\,v_0\,, \label{iged1}\\
z(t) &=& t\,z_0\,, \label{iged2}\\
u(t) &=& t\,u_0 + \half t^2 \,\sJ v_0\,.\eop \label{iged3}
\end{eqnarray}
\end{corollary}
\begin{corollary}
If\label{dgc} $[\n,\n] \subseteq \U$ and $\E=\{0\}$, then $N$ is
geodesically connected.  Consequently, so is any nilmanifold with such a
universal covering space.
\end{corollary}
\begin{proof}
$N$ is complete by Theorem \ref{gt}, hence pseudoconvex.  The preceding
geodesic equations show that $N$ is nonreturning.  Thus the space of
geodesics $G(N)$ is Hausdorff by Theorem 5.2 of \cite{BP9}.  Now Theorem
4.2 of \cite{BLP} yields geodesic connectedness of $N$.
\end{proof}
Thus these compact nilmanifolds are much like tori.  This is also
illustrated by the computation of their period spectrum in Theorem
\ref{upd}.

When $J^2$ is diagonalizable,                  %when *does* this happen?
then the nondegenerate case also can be made completely explicit.  For
example, this happens when $\langle\,,\rangle$ restricted to $\E_2$ is
definite ({\em e.\,g.,} in the Riemannian case).  The extreme instance is
when $J^2$ is a scalar multiple of the identity.  In the Riemannian case,
this is much the same as $H\!$-type; however, in general this condition is
not related to \ph-type.

Let $\{\th_1,\ldots,\th_k\}$ be the set of distinct nonzero
eigenvalues of $J^2$. Decompose $\E_2$ as the orthogonal direct sum
$\bigoplus_{j=1}^k \w_j$ where $J$ leaves each $\w_j$ invariant and
$J^2|_{\w_j}=\th_j I$.  Write $e_0=e_1+e_2$ with $e_i\in\E_i$ and then
write $e_2=\sum_{j=1}^k w_j$ with each $w_j\in\w_j$.
\begin{proposition}
If\label{igendd} $N$ has a nondegenerate center and $J^2$ diagonalizes, then
\begin{eqnarray}
e(t) &=& t\,e_1 + \left( e^{tJ} - I\right) J^{-1}e_2\,,\label{igendd1}\\
z(t) &=& t\,z_1(t) + z_2(t)\,,\label{igendd2}
\end{eqnarray}
where
\begin{eqnarray*}
z_1(t) &=& z_0+\half\left[e_1,\left(e^{tJ}+I\right)J^{-1}e_2\right]
           +\half\sum_{j=1}^k\left[J^{-1}w_j,w_j\right],\\
z_2(t) &=& \left[e_1,\left(I-e^{tJ}\right)\left(
      J^{-2}e_2\right)\right]+\half\left[e^{tJ}J^{-1}e_2,
      J^{-1}e_2\right]\\
&& {}-\half\sum_{i\ne j}\frac{1}{\th_j-\th_i}\Bigl(\left[e^{t
   J}Jw_i,e^{tJ} J^{-1}w_j\right]-\left[e^{tJ}w_i,
   e^{tJ}w_j\right]\Bigr)\\
&& {}+\half\sum_{i\ne j}\frac{1}{\th_j-\th_i}\Bigl(\left[Jw_i,
J^{-1}w_j\right]-\left[w_i,w_j\right]\Bigr).
\end{eqnarray*}
\end{proposition}
\begin{proof}
Compare Eberlein \cite[(3.5)]{E'}. First note that on each $\w_j$ where the
restriction $J$ is invertible, $J = \th_j J^{-1}$.

Next, one verifies straightforwardly that
\begin{equation}
[\dot{e}(t),e(t)] = \left[e_1,\left(e^{tJ}-I\right)J^{-1}
e_2\right] + t\left[e^{tJ}e_2,e_1\right] + \left[e^{tJ}e_2,\left(
e^{tJ}-I\right)J^{-1} e_2\right],
\label{e-c}
\end{equation}
and, using the properties of $J$, that the derivative with
respect to $t$ of each $\left[e^{tJ}w_j,e^{tJ}J^{-1}w_j\right]$ is zero.
Hence
$$
\left[e^{tJ}w_j,e^{tJ}J^{-1}w_j\right]=\left[w_j,J^{-1}w_j\right]
$$
for all $t$ and $j$. From this and the decomposition of $e_2$, we obtain
\begin{equation}
-2\dot{z}_3 = \left[e^{tJ}e_2,e^{tJ}J^{-1}e_2\right] =
\sum_{i\ne j}\left[e^{tJ}w_i,e^{tJ}J^{-1}w_j\right] +
\sum_j\left[w_j,J^{-1}w_j\right]. \label{e-e}
\end{equation}
Using the two facts about $J$ again, another computation produces
$$ \dot{z}_2 = -\left[e_1,e^{tJ}J^{-1}e_2\right] +
\half\left[e^{tJ}e_2,J^{-1}e_2\right]. $$
Combining this with (\ref{e-e}) and $\dot{z}(t)=z_1(t) + t\,\dot{z}_1 +
\dot{z}_2 + \dot{z}_3$, we get
\begin{eqnarray*}
\dot{z}(t) &=& z_0 - \half\left[e^{tJ}e_2,e^{tJ}J^{-1}e_2\right] -
\half\left[e_1,\left(e^{tJ}-I\right)J^{-1}e_2\right] \\
&& {}+\half\left[e^{tJ}e_2,J^{-1}e_2\right] -
{\txs\frac{t}{2}}\left[e^{tJ}e_2,e_1\right].
\end{eqnarray*}
Finally, it follows from (\ref{e-c}) and this that $z$ satisfies
$\dot{z} + \half[\dot{e},e] = z_0$, and from (\ref{igendd2}) and this that
$z$ satisfies the initial conditions.
\end{proof}
Next we recover \cite[(3.7)]{E'} in our case.
\begin{corollary}
Under\label{csgf} the same hypotheses,  
$e(t)=t\,e_1 + \sum_j\xi_j(t)$ where
$$
\xi_j(t) = \left\{ \begin{array}{lcl}
 \dps\left(\cos t\l_j - 1\right)J^{-1}w_j + \frac{\sin
 t\l_j}{\l_j}w_j & & \dps\mbox{for \,} \th_j=-\l_j^2\,,\\[3ex]
 \dps\left(\cosh t\l_j - 1\right)J^{-1}w_j + \frac{\sinh
 t\l_j}{\l_j}w_j & \quad & \dps\mbox{for \,} \th_j=\l_j^2\,,
\end{array}\right.
$$
with $\l_j>0$.
\end{corollary}
\begin{proof}
First use the decomposition of $e_2$. From the representation of $J$ on 
$\w_j$, we obtain in the first case
$$
e^{tJ} = \cos t\l_jI + \frac{\sin t\l_j}{\l_j}J
$$
there, and the second case is similar.
\end{proof}

Thus some geodesics in the pseudoriemannian case may exhibit the same
helical behavior as geodesics do in the Riemannian case.  Different
signatures on the Heisenberg group illustrate both possibilities (for
example, $+--$ and $--+$, respectively).
\begin{remark}
If\label{sgf} $a_0 = 0$ or $x_0 = 0$, then the geodesic becomes $\g(t) =
\exp(tx_0)$ or $\exp(ta_0)$, respectively.  More generally, it follows
from the geodesic equation (\ref{ges}) that if $\g(t)$ is the unique
geodesic with $\dg(0) = a_0 + x_0$, then $\g(t) = \exp\left( t\left( a_0 +
x_0 \right)\right)$ if and only if $\add{x_0}{(a_0 + x_0)} = 0$ if and only
if $(a_0 + x_0)\perp [x_0,\n]$.  We note that this is also equivalent to
$\add{x_0}{(z_0 + v_0)} = 0$ if and only if $(z_0 + v_0) \perp [x_0,\n]$.
\end{remark}
\begin{proposition}
Let $N$ be a simply connected, 2-step nilpotent Lie group with
left-invariant metric tensor $\langle\,,\rangle$, and let $\g$ be a
geodesic with $\dg(0) = a\in \z$. If $\add{\bul}{a} = 0$, then there are
no conjugate points along $\g$.
\end{proposition}
\begin{proof}
It follows from Theorems \ref{zph} and \ref{kze} and Proposition
\ref{oscn} that all planes containing any $\dg(t)$ are homaloidal.  Thus
the Jacobi equation trivializes and the result follows.
\end{proof}
In the Riemannian case \cite[(3.10)]{E'} and in the timelike Lorentzian
case (using \cite[p.278, 23.\,Lemma]{O} and an argument {\em mutatis
mutandis\/}), the converse holds.  It is not clear if it can be
generalized much more than this.

\section{Isometry group}
\label{isogp}

We begin with a general result about nilpotent Lie groups.
\begin{lemma}
Let\label{iaid} $N$ be a connected, nilpotent Lie group with a left-invariant 
metric tensor. Any isometry of $N$ that is also an inner automorphism must 
be the identity map.
\end{lemma}
\begin{proof}
Consider an isometry of $N$ which is also the inner automorphism
determined by $g\in N$. Then $\Ad{g}$ is a linear isometry of \n. The Lie
group exponential map is surjective, so there exists $x\in\n$ with
$g=\exp(x)$. Then
$$ \Ad{g}=e^{\ad{x}}. $$
Now $\ad{x}$ is nilpotent so all its eigenvalues are 0. Thus all the
eigenvalues of $\Ad{g}$ are 1, so it is unipotent. It now follows
that $\Ad{g}$ is the identity. Thus $g$ lies in the center of $N$ and
the isometry of $N$ is the identity.
\end{proof}

Now we give some general information about the isometries of $N$. Letting
$\Aut(N)$ denote the automorphism group of $N$ and $I(N)$ the isometry
group of $N$, set $O(N) = \Aut(N)\cap I(N)$. In the Riemannian case, $I(N)
= O(N)\ltimes N$, the semidirect product where $N$ acts as left
translations.  We have chosen the notation $O(N)$ to suggest an analogy
with the pseudoeuclidean case in which this subgroup is precisely the
(general, including reflections) pseudorthogonal group.  According to
Wilson \cite{Wi}, this analogy is good for any nilmanifold (not necessarily
2-step).

To see what is true about the isometry group in general, first consider the
(left-invariant) splitting of the tangent bundle $TN = \z N \ds \v N$.
\begin{definition}
Denote\label{displ} by $\Ispl(N)$ the subgroup of the isometry group
$I(N)$ which preserves the splitting $TN = \z N \ds \v N$. Further, let 
$\Iaut(N) = O(N)\ltimes N$, where $N$ acts by left translations.
\end{definition}
\begin{proposition}
If\label{sispl} $N$ is a simply-connected, 2-step nilpotent Lie group with
left-invariant metric tensor, then $\Ispl(N) \le \Iaut(N)$.
\end{proposition}
\begin{proof}
Similarly to the observation by Eberlein \cite{E'}, it is now easy to check
that the relevant part of Kaplan's proof for Riemannian $N$ of $H\!$-type
\cite{K} is readily adapted and extended to our setting, with the proviso
that it is easier to just replace his expressions $j(a)x$ with $\add{x}a$
instead of trying to convert to our $j$.
\end{proof}
We shall give an example to show that $\Ispl < \Iaut$ is possible when $\U 
\ne \{0\}$.

When the center is degenerate, the relevant group analogous to a
pseudorthogonal group may be larger.
\begin{proposition}
Let\label{tO} $\wO(N)$ denote the subgroup of $I(N)$ which fixes $1 \in N$.
Then $I(N) \cong \wO(N) \ltimes N$, where $N$ acts by left translations.\eop
\end{proposition}
The proof is obvious from the definition of $\wO$.  It is also obvious that 
$O \le \wO$.  We shall give an example to show that $O < \wO$, hence $\Iaut
< I$, is possible when the center is degenerate.

Thus we have three groups of isometries, not necessarily equal in general:
$\Ispl \le \Iaut \le I$.  When the center is nondegenerate ($\U = \{0\}$),
the Ricci transformation is block-diagonalizable and the rest of Kaplan's
proof using it now also works.
\begin{corollary}
If\label{isplnd} the center is nondegenerate, then $I(N) = \Ispl(N)$ 
whence $\wO(N) \cong O(N)$.\eop
\end{corollary}

We now begin to prepare for the promised examples.
Recall the result in \cite[3.57(b)]{Wa} that for any simply connected Lie
group $G$ with Lie algebra $\frak{g}$,
$$ \Aut(G) \cong \Aut(\frak{g}) : \alpha\longmapsto\alpha_*|_{\frak{g}}\,.$$
Now, if $\alpha$ is also an isometry, then so is $\alpha_*$. Conversely, any
isometric automorphism of $\frak{g}$ comes from an isometric automorphism of
$G$ by means of left-invariance (or homogeneity). Therefore
$$ O(N) = \Aut(N) \cap I(N) \cong \Aut(\n) \cap O^p_q $$
for one of our simply connected groups $N$ with signature $(p,q)$. So up to
isomorphism, we may regard
$$ O(N) \le O^p_q \,.$$
\begin{example}
Any\label{sf} flat, simply connected $N$ is a spaceform.  (By Corollary
\ref{ccf} no other constant curvature is possible.)  Thus the isometry group
is known (up to isomorphism) and $\wO(N) \cong O^p_q$ for any such $N$ of
signature $(p,q)$.
\end{example}

For the rest of the examples, we need some additional preparations.  Let
$\h_3$ denote the 3-dimensional Heisenberg algebra and consider the unique
\cite{D}, 4-dimensional, 2-step nilpotent Lie algebra $\n_4 = \h_3\times\R$,
with basis $\{v,e,z,u\}$, structure equation $[v,e] = z$, and center $\z =
\lsp z,u\rsp$.

A general automorphism of $\n_4$ as a vector space, $f \colon \n_4 \to \n_4$, will
be given with respect to the basis $\{v,e,z,u\}$ by a matrix
$$ \left[
\begin{array}{cccc}
a_1 & b_1 & c_1 & d_1\\
a_2 & b_2 & c_2 & d_2\\
a_3 & b_3 & c_3 & d_3\\
a_4 & b_4 & c_4 & d_4
\end{array} \right] .$$
In order for $f$ to be a Lie algebra automorphism of $\n_4$, it must preserve
the center, and $[f(v),f(e)]=f(z)$ while the other products of images by $f$
are zero; all this together implies that the matrix of $f$ must have the
form
$$ \left[
\begin{array}{cccc}
a_1 & b_1 &  0  & 0 \\
a_2 & b_2 &  0  &  0 \\
a_3 & b_3 & c_3 & d_3\\
a_4 & b_4 &  0  & d_4
\end{array} \right] $$
with $c_3=a_1\,b_2 - a_2\,b_1\not=0$, $d_4\not=0$ and arbitrary
$a_3, a_4, b_3, b_4, d_3$.

Now the matrix of $f$ can be decomposed as a product of matrices $M_1 \cdot
M_2 \cdot M_3 \cdot M_4$, where
$$ \begin{array}{lcl}
  M_1 = \left[
   \begin{array}{cccc}
     a_1 & b_1 &  0  &  0 \\
     a_2 & b_2 &  0  &  0 \\
      A  &  B  & c_3 &  0 \\
      0  &  0  &  0  &  1
   \end{array} \right], & \quad &
  M_2 = \left[
   \begin{array}{cccc}
     1 & 0 &  0  &  0 \\
     0 & 1 &  0  &  0 \\
     0 & 0 &  1  &  0 \\
     0 & 0 &  0  & d_4
   \end{array} \right],\\[7ex]
  M_3 = \left[
   \begin{array}{cccc}
      1  &  0  &  0  &  0 \\
      0  &  1  &  0  &  0 \\
      0  &  0  &  1  &  d_3/c_3 \\
      0  &  0  &  0  &  1
   \end{array} \right], & &
  M_4 = \left[
   \begin{array}{cccc}
     1 & 0 &  0  &  0 \\
     0 & 1 &  0  &  0 \\
     0 & 0 &  1  &  0 \\
     C & D &  0  &  1
   \end{array} \right],
\end{array} $$
with
$$ A=\frac{a_3\,d_4 - a_4\,d_3}{d_4}\,,\ \
B=\frac{b_3\,d_4 - b_4\,d_3}{d_4}\,,\ \
C=\frac{a_4}{d_4}\,,\ \ D=\frac{b_4}{d_4}\,, $$
which are well defined because $c_3,d_4\not= 0$.
In this decomposition of $f$ we have:
$$ \begin{array}{l}
M_1 \in \Aut(\h_3),\ \mbox{and they fill out } \Aut(\h_3);\\[1ex]
M_2 \in \Aut(\R);\\[1ex]
M_3 \in M \ \mbox{mixes part of the center but fixes $v,e$, and $z$;}\\[1ex]
M_4 \ \mbox{fixes $z$ and $u$ but the image of $\lsp v,e\rsp$ is contained
   in $\lsp v,e,u\rsp$.}
\end{array} $$
The matrices $M_4$ form a group isomorphic to an abelian $\R^2$. Note that
for such an $M_4$, $v\mapsto v + Cu$ and $e\mapsto e + Du$.
\begin{proposition}
The decomposition of $\Aut(\n_4)$ is
$$ \Aut(\n_4) = \Aut(\h_3)\cdot \Aut(\R)\cdot M\cdot \R^2 .\eop$$
\end{proposition}

Now we compute $O(N)$ for several simply-connected groups of interest.
Specifically, we shall compute $O(N) \cong O(\n)$ for the flat and
nondegenerate 3-dimensional Heisenberg groups $H_3$ and the 4-dimensional
groups $N_4 \cong H_3\times \R$ with all possible signatures.

\begin{example}
We begin with the form for automorphisms of $\h_3$ on the basis $\{v,e,u\}$,
$$ f = \left[\begin{array}{ccc}
   a_1 & b_1 &  0 \\
   a_2 & b_2 &  0 \\
   a_3 & b_3 & c_3 \end{array} \right] $$
with $c_3=a_1\,b_2 - a_2\,b_1 \neq 0$ and arbitrary $a_3,b_3$.  For the flat
metric tensor we use
$$ \eta = \left[\begin{array}{ccc}
   0 & 0 & 1 \\
   0 & \bve & 0 \\
   1 & 0 & 0 \end{array} \right] $$
and we need to solve
$$ f^T\eta f = \eta\,. $$
Using {\it Mathematica\/} to help, it is easy to see that we must have $b_1
= 0$, whence $a_1 \neq 0$ and thus the \verb+Solve+ command in fact does
give all the solutions.  We also obtain $a_1 = \pm 1$, so
\begin{equation}
\arraycolsep.5em
f = \left[ \begin{array}{ccc}
   \pm 1 & 0 & 0 \\
   a_2 & 1 & 0 \\
   \mp\bve a_2^2/2 & \mp\bve a_2 & \pm 1 \end{array} \right]
\label{d3d}
\end{equation}
with arbitrary $a_2$.  We note that the determinant is 1, so that this
1-parameter group lies in $SO^1_2$ or $SO^2_1$ according as $\bve = \pm 1$,
respectively.  With some additional work, this is seen to be a group of
horocyclic translations ($HT$ in \cite{CP1}) subjected to a change of basis
(rotation and reordering).  (Comparing parameters, we have here $a_2 =
-t\sqrt{2}$ there.)  Note $\dim O(H_3) = 1$ while from Example \ref{sf} we
have $\dim\wO(H_3) = 3$.  We can identify the rest of $\wO$ from an Iwasawa
decomposition:  $O$ is the nilpotent part, so the rest of $\wO$ consists of
rotations and boosts.
\end{example}

Many groups with degenerate center have such a subgroup isometrically
embedded.
\begin{proposition}
Assume\label{ifh} $\U\ne\{0\}\ne\E$ and let $u\in\U$, $v\in\V$, such that
$\langle u,v\rangle = 1$ and $0 \ne j(u)v = e \in \E$.  Then $j(u)^2 = -I$
on $\lsp v,e\rsp$ if and only if $\h = \lsp u,v,e\rsp$ is isometric to the
3-dimensional Heisenberg algebra with null center. Consequently, $H =
\exp\h$ is an isometrically embedded, flat subgroup of $N$.
\end{proposition}
\begin{proof}
It is easy to see that $j(u)e = -v$; {\em cf.} Remarks \ref{rems}. Then,
on the one hand, $\langle[v,j(u)v],v\rangle = \langle[v,e],v\rangle$.  On
the other hand, $\langle[v,j(u)v],v\rangle = -\langle[j(u)v,v],v\rangle =
-\langle\add{j(u)v}{\io u},v\rangle = -\langle\io j(u)^2 v,v\rangle =
\langle\io v,v\rangle = \langle u,v\rangle = 1$. It follows that $[v,e] =
u$.  The converse follows from Example \ref{h3n}.
\end{proof}
\begin{corollary}
For\label{ifhi} any $N$ containing such a subgroup,
$$\Ispl(N) < \Iaut(N) < I(N) .\eop$$
\end{corollary}
Unfortunately, this class does not include our flat groups of Theorem
\ref{e0f} in which $[\n,\n] \subseteq \U$ and $\E = \{0\}$.  However, it
does include many groups that do not satisfy $[\n,\n] \subseteq \U$, such
as the simplest quaternionic Heisenberg group of Example \ref{hq} {\em et
seq.}
\begin{remark}
A\label{kgt} direct computation shows that on this flat $H_3$ with null
center, the only Killing fields with geodesic integral curves are the
nonzero scalar multiples of $u$.
\end{remark}

Similarly, we can do the 3-dimensional Heisenberg group with nondegenerate
center (including certain definite cases).  With respect to the orthonormal
basis $\{e_1,e_2,z\}$ with structure equation $[e_1,e_2] = z$, we obtain
$$ \wO(H_3) = O(H_3) \cong \left\{ \begin{array}{rl}
      O_2\quad & \mbox{if } \bve_1\bve_2 = 1,\\[1ex]
      O^1_1\quad & \mbox{if } \bve_1\bve_2 = -1 \end{array}\right. $$
of dimension two.

In general, we always can arrange to have $e_1$ and $e_2$ orthonormal, but
we may not be able to keep the structure equation {\em and\/} have $z$ a
unit vector simultaneously.  In the indefinite case, we always can have the
$z$-axis orthogonal to the $e_1e_2$-plane.  But in the definite case, the
$e_1e_2$-plane need not be orthogonal to the $z$-axis.  It is not yet clear
how best to compute $O(H_3)$ in these cases.

\begin{example}
Now we consider the flat $N_4 \cong H_3 \times \R$ with null center.
Here we have the basis $\{v_1,v_2,u_1,u_2\}$ of $\n_4$ with structure 
equation $[v_1,v_2] = u_1$, the form for automorphisms
$$ f = \left[ \begin{array}{cccc}
   a_1 & b_1 &  0  & 0 \\
   a_2 & b_2 &  0  &  0 \\
   a_3 & b_3 & c_3 & d_3\\
   a_4 & b_4 &  0  & d_4 \end{array} \right] $$
with $c_3=a_1\,b_2 - a_2\,b_1\not=0$, $d_4\not=0$, and arbitrary
$a_3,a_4,b_3,b_4,d_3$, and the flat metric tensor
$$ \eta = \left[\begin{array}{cccc}
   0 & 0 & 1 & 0 \\
   0 & 0 & 0 & 1 \\
   1 & 0 & 0 & 0 \\
   0 & 1 & 0 & 0 \end{array}\right] .$$
Again we find $b_1 = 0$ and $a_1 \neq 0$ and \verb+Solve+ gives all
solutions of $f^T\eta f = \eta$. We obtain
\begin{equation}
\arraycolsep.5em
   f = \left[\begin{array}{cccc}
   a_1 & 0 & 0 & 0 \\
   a_2 & 1/a_1^2 & 0 & 0 \\
   a_1^2 a_2 b_3 & b_3 & 1/a_1 & -a_1 a_2 \\
   -a_1^3 b_3 & 0 & 0 & a_1^2 \end{array} \right]
\label{d422}
\end{equation}
with $a_1\neq 0$ and arbitrary $a_2,b_3$. Again we note the determinant is 1
so this 3-parameter group lies in $SO^2_2$.  Note $\dim O(N_4) = 3$ while
from Example \ref{sf} we have $\dim\wO(N_4) = 6$.

We consider three 1-parameter subgroups.
$$ \begin{array}{lcl}
   a_2 = b_3 = 0 & \qquad & \mbox{I} \\
   a_1 = 1,\; b_3 = 0 && \mbox{II} \\
   a_1 = 1,\; a_2 = 0 && \mbox{III} \end{array} $$
We make a $2+2$ decomposition of our space into $\lsp v_1,u_1\rsp \ds \lsp
v_2,u_2\rsp$.  If we think of $a_1 = \pm e^t$, then subgroup I is a boost by
$t$, respectively $2t$, on the two subspaces.  Subgroups II and III are both
horocyclic translations on the two subspaces, and they commute with each
other.

Again, we can see what the part of $\wO$ outside $O$ looks like from an
Iwasawa decomposition $\wO = K\!AN \cong O^2_2$.  Clearly, I is a subgroup
of $A$ while II and III are subgroups of $N$.  Now $O^2_2$ has rank 2 so
$\dim A = 2$.  Since $\dim K = 2$, then $\dim N = 2$ also.  Thus all the
rotations and half of the boosts are the part of $\wO$ outside $O$.
\end{example}

Many of our flat groups from Theorem \ref{e0f} have such a
subgroup isometrically embedded, as in fact do many others which
are not flat.
\begin{proposition}
Assume\label{if4} $[\V,\V] \subseteq \U$ and $\dim\U = \dim\V \ge 2$.  Let
$u_1,u_2 \in \U$ and $v_1,v_2 \in \V$ be such that $\langle u_i,v_i\rangle =
1$.  If, on $\lsp v_1,v_2\rsp$, $j(u_1)^2 = -I$ and $j(u) = 0$ for $u \notin
\lsp u_1\rsp$, then $\lsp u_1,u_2,v_1,v_2\rsp \cong \h_3\times\R$.
\end{proposition}
\begin{proof}
We may assume without loss of generality that $j(u_1)v_1 = v_2$.
Then a straightforward calculation shows that $[v_1,v_2] = u_1$
and the result follows.
\end{proof}
\begin{corollary}
For\label{if4i} any $N$ containing such a subgroup,
$$\Ispl(N) < \Iaut(N) < I(N) .\eop$$
\end{corollary}

\begin{example}
Next we consider some $N_4$ with only partially degenerate center.  Here we
have the basis $\{v,e,z,u\}$ of $\n_4$ with structure equation $[v,e] = z$
and signature $(+\,\bve\,\ve\,-)$.  We take $f$ again as above, and the
metric tensor
$$ \eta = \left[\begin{array}{cccc}
   0 & 0 & 0 & 1 \\
   0 & \bve & 0 & 0 \\
   0 & 0 & \ve & 0 \\
   1 & 0 & 0 & 0 \end{array}\right] .$$
Applying \verb+Solve+ to $f^T\eta f = \eta$ directly is quite wasteful this
time, producing many solutions with $\ve = 0$ which must be discarded, and
producing duplicates of the 4 correct solutions.  After removing the
rubbish, however, it is easy to check that there are no more again.  Thus we
obtain
\begin{equation}
\arraycolsep .5em
\begin{array}{c}
f = \left[\begin{array}{cccc}
   \pm 1 & 0 & 0 & 0 \\
   a_2 & \pm 1 & 0 & 0 \\
   0 & 0 & 1 & 0 \\
   \mp\bve a_2^2/2 & \mp\bve a_2 & 0 & \pm 1 \end{array} \right]\\[7ex]
\mbox{or}\\[2ex]
f = \left[\begin{array}{cccc}
   \pm 1 & 0 & 0 & 0 \\
   a_2 & \mp 1 & 0 & 0 \\
   0 & 0 & -1 & 0 \\
   \mp\bve a_2^2/2 & \mp\bve a_2 & 0 & \pm 1 \end{array} \right]
\end{array}
\label{d4os}
\end{equation}
with arbitrary $a_2$.  Now we have the determinant $\pm 1$ for these
1-parameter subgroups, and both are contained in $O^{1+}_3$, $O^{2+}_2$, or
$O^{3+}_1$, according to the choices for $\bve$ and $\ve$, respectively
$\bve=\ve=1$, $\bve\neq\ve$, or $\bve=\ve= -1$.  Again, these are all
horocyclic translations as found in \cite{CP1}, modulo change of basis.

This yields $\dim O(N_4) = 1$, while in the next example we shall see that 
$\dim \wO(N_4) = 2$. % embedding result a la {ifh} and {if4}? 
\end{example}

Since $N$ is complete, so are all Killing fields.  It follows
\cite[p.\,255]{O} that the set of all Killing fields on $N$ is the Lie
algebra of the (full) isometry group $I(N)$.  Also, the 1-parameter groups
of isometries constituting the flow of any Killing field are global. Thus
integration of a Killing field produces (global, not just local) isometries
of $N$.
\begin{example}
Consider again the unique \cite{D}, nonabelian, 2-step nilpotent Lie algebra
$\n_4$ of dimension 4. We take a metric tensor so that the center is only
partially degenerate, as in the preceding Example.  Thus we choose the basis
$\{v,e,z,u\}$ with structure equation $[v,e] = z$ and nontrivial inner
products
$$ \langle v,u\rangle = 1,\qquad\langle
e,e\rangle = \bve,\qquad\langle z,z\rangle = \ve.$$
This basis is partly null and partly orthonormal.  The associated simply
connected group $N_4$ is isomorphic to $H_3\times\R$ and we realize it as
the real matrices
\begin{equation}
\left[\begin{array}{ccccc}
      1 & x & z & 0 & 0 \\
      0 & 1 & y & 0 & 0 \\
      0 & 0 & 1 & 0 & 0 \\
      0 & 0 & 0 & 1 & w \\
      0 & 0 & 0 & 0 & 1 \end{array}\right]
\end{equation}
so this group is diffeomorphic to $\R^4$.  Take
\begin{equation}
v = \frac{\partial}{\partial x}\,,\quad
e = \frac{\partial}{\partial y} + x\frac{\partial}{\partial z}\,,\quad
z = \frac{\partial}{\partial z}\,,\quad
u = \frac{\partial}{\partial w}\,,\quad
\end{equation}
which realizes the structure equation in the given representation. (The
double use of $z$ should not cause any confusion in context.)   The
nontrivial adjoint maps are
$$ \add{v}{z} = \ve\bve\,e \qquad\mbox{and}\qquad \add{e}{z} = -\ve\,u $$
so the nontrivial covariant derivatives are
\begin{eqnarray*}
\del{v}e = -\del{e}v &=& \half z\,,\\
\del{z}v = \del{v}z &=& -\half\ve\bve\,e\,,\\
\del{z}e = \del{e}z &=& \half\ve\,u\,.
\end{eqnarray*}
\begin{remark}
The space $N_4$ is not flat; in fact, a direct computation shows that
$$ \langle R(z,v)v,z \rangle = \quar\bve\qquad\mbox{and}\qquad
\langle R(v,e)e,v \rangle = -{\mbox{$\txs\frac{3}{4}$}}\ve\, .$$
\end{remark}

We consider Killing fields $X$ and write $X = f^1 u + f^2 z + f^3 v + f^4
e$.  We denote partial derivatives by subscripts.  From Killing's equation
we obtain the system.
\begin{eqnarray}
f^1_x = f^2_z = f^3_w &=& 0 \label{kf1}\\
f^1_w + f^3_x &=& 0 \label{kf2}\\
\ve f^2_w + f^3_z &=& 0 \label{kf3}\\
f^4_y + x f^4_z &=& 0 \label{kf4}\\
\ve f^1_z + f^2_x + f^4 &=& 0 \label{kf5}\\
\bve f^4_x + f^1_y + x f^1_z &=& 0 \label{kf6}\\
\bve f^4_w + f^3_y +x f^3_z &=& 0 \label{kf7}\\
\ve\bve f^4_z + f^2_y &=& f^3 \label{kf8}
\end{eqnarray}
Simple considerations of second derivatives yield these relations.
\begin{eqnarray*}
f^1_x = f^1_{ww} = f^1_{wy} = f^1_{wz} &=& 0 \\
f^2_z = f^2_{ww} = f^2_{wx} &=& 0 \\
f^3_w = f^3_{xx} = f^3_{xz} = f^3_{zz} &=& 0 \\
f^4_w = f^4_{xz} &=& 0
\end{eqnarray*}
From the first we get $f^1_w = A$, and then from (\ref{kf2}) $f^3_x = -A$,
so in particular $f^3_{xy}=0$.  Then (\ref{kf7}) and $f^4_w = f^3_{xy} =
f^3_{xz} = 0$ imply $f^3_z = 0$ and therefore $f^3_y = 0$ also.  Hence,
$$f^3 = -Ax + C_1\, .$$

Also, (\ref{kf3}) and $f^3_z = 0$ imply $f^2_w =0$, and (\ref{kf8}) and
$f^4_{zx} = 0$ imply $f^2_{yx} = f^3_x = -A$.  It then follows that
\begin{eqnarray*}
f^1 &=& A w + F^1(y,z) \, ,\\
f^2 &=& F^2(x,y) \, , \quad \mbox{with } F^2_{xy} = -A \, ,\\
f^3 &=& -A x + C_1 \, ,\\
f^4 &=& f^4(x,y,z) \, , \quad \mbox{with } f^4_{xz} = 0 \, .
\end{eqnarray*}
Now, using (\ref{kf6}) and $f^4_{xz}=0$, after a few steps we obtain
$F^1_z = C_2$ whence
$$F^1(y,z) = C_2 z + G^1 (y) \, .$$
Then (\ref{kf5}) gives us $f^4_z = - \ve f^1_{zz} = 0$ so
$$f^4 = - \ve C_2 - F^2_x (x,y) \, .$$

From (\ref{kf4}) and $f^4_z = 0$ we get $0=f^4_y=-F^2_{xy}=A$ whence
$$ A=0 \, .$$
From (\ref{kf8}) and $f^4_z=0$ we get $f^2_y = F^2_y = C_1$ whence
$F^2(x,y) = C_1 y + G^2(x)$;
hence
\begin{eqnarray*}
f^1 &=& C_2 z + G^1(y) \, ,\\
f^2 &=& C_1 y + G^2(x) \, ,\\
f^3 &=& C_1 \, ,\\
f^4 &=& -\ve C_2 - G^2_x (x)  \, .
\end{eqnarray*}

Finally,
(\ref{kf6}) implies $-\bve G^2_{xx} + G^1_y +C_2 x = 0$ so
$-\bve G^2_{xxx} = -C_2$
and hence
\begin{eqnarray*}
G^2_{xxx} &=& \bve C_2 \, ,\\
G^2_{xx} &=& \bve C_2 x + C_3 \, ,\\
G^2_x &=& \frac{\bve C_2}{2} x^2 + C_3 x + C_4 \, ,\\
G^2 &=& \frac{\bve C_2}{6} x^3 + \frac{C_3}{2} x^2 + C_4 x + C_5 \, .
\end{eqnarray*}
Then $G^1_y = \bve G^2_{xx} - C_2 x = \bve C_3$ whence
$$G^1(y) = \bve C_3 y + C_6 \, .$$
and we are done.
\begin{eqnarray*}
f^1 &=& \bve C_3 y + C_2 z + C_6 \\
f^2 &=& \frac{\bve C_2}{6} x^3 + \frac{C_3}{2} x^2 + C_4 x + C_1 y + C_5 \\
f^3 &=& C_1 \\
f^4 &=& -\frac{\bve C_2}{2} x^2 - C_3 x - (\ve C_2 + C_4)
\end{eqnarray*}
Therefore, $\dim I(N_4) = 6$.

The Killing vector fields with $X(0)=0$ arise for $C_1 = C_5 = C_6 = 0$ and
$C_4 = -\ve C_2$, so their coefficients are
\begin{eqnarray*}
f^1 &=& \bve C_3 y + C_2 z \, ,\\
f^2 &=& \frac{\bve C_2}{6} x^3 + \frac{C_3}{2} x^2 - \ve C_2 x \, ,\\
f^3 &=& 0 \, ,\\
f^4 &=& -\frac{\bve C_2}{2} x^2 - C_3 x \, .
\end{eqnarray*}
Therefore $\dim\wO(N_4) = 2$.
\begin{proposition}
The Killing vector fields $X$ with $X(0) = 0$ are
\begin{eqnarray*}
X &=& -\left( \frac{\bve C_2}{2} x^2 + C_3 x \right)
   \frac{\partial}{\partial y} - \left( \frac{\bve C_2}{3} x^3 +
   \frac{C_3}{2} x^2 + \ve C_2 x \right)\frac{\partial}{\partial z}\\
&&{}+ \left(\bve C_3 y + C_2 z \right)\frac{\partial}{\partial w} \, .
\end{eqnarray*}
\end{proposition}

Their integral curves are the solutions of
\begin{eqnarray*}
\dot{x}(t) &=& 0 \, ,\\
\dot{y}(t) &=& -\frac{\bve C_2}{2} x^2 - C_3 x \, ,\\
\dot{z}(t) &=& -\frac{\bve C_2}{3} x^3 - \frac{C_3}{2} x^2 - \ve C_2 x \,
   ,\\
\dot{w}(t) &=& \bve C_3 y + C_2 z \, ,
\end{eqnarray*}
and the solutions with initial condition $(x_0,y_0,z_0,w_0)$ are
\begin{eqnarray*}
x(t) &=& x_0 \, ,\\
y(t) &=& -\left( \frac{\bve C_2}{2} x^2_0 + C_3 x_0 \right) t + y_0 \, ,\\
z(t) &=& -\left( \frac{\bve C_2}{3} x^3_0 + \frac{C_3}{2} x^2_0 + \ve C_2
   x_0 \right) t + z_0 \, ,\\
w(t) &=& -\left( \frac{\bve C_2^2}{3} x^3_0 + C_2 C_3 x^2_0 + (\ve C_2^2 +
   \bve C_3^2)x_0 \right) \frac{t^2}{2} + (\bve C_3 y_0 + C_2 z_0) t + w_0
   \, .
\end{eqnarray*}

Thus the 1-parameter group $\F$ of isometries is given by
\begin{eqnarray*}
\F^1 &=& x \, ,\\
\F^2 &=& -\left( \frac{\bve C_2}{2} x^2 + C_3 x \right) t + y \, ,\\
\F^3 &=& -\left( \frac{\bve C_2}{3} x^3 + \frac{C_3}{2} x^2 + \ve C_2 x
   \right) t + z \, ,\\
\F^4 &=& -\left( \frac{\bve C_2^2}{3} x^3 + C_2 C_3 x^2 + (\ve C_2^2 + \bve
   C_3^2)x \right)\frac{t^2}{2} + (\bve C_3 y + C_2 z) t + w \, ,
\end{eqnarray*}
and its Jacobian is
$$ \arraycolsep .7em
\F_* = \left[ \begin{array}{cccc}
1 & 0 & 0 & 0\\
-(\bve C_2x + C_3)t & 1 & 0 & 0\\
-(\bve C_2 x^2 + C_3x +\ve C_2)t & 0 & 1 & 0\\
-\half(\bve C_2^2x^2 +2C_2C_3x+(\ve C_2^2 +\bve C_3^2)) t^2 & \bve C_3 t &
   C_2t & 1
 \end{array}\right] .$$
Hence
\begin{eqnarray*}
\F_* v &=& v - (C_3 + \bve C_2 x) t\,e - \ve C_2 t \, z - (\bve C_2^2
   x^2 + 2 C_2C_3 x + \ve C_2^2 + \bve C_3^2 ) \frac{t^2}{2}\, u\, ,\\
\F_* e &=& e + (\bve C_3 + C_2 x)t \, u\, ,\\
\F_* z &=& z + C_2 t\,u\,,\\
\F_* u &=& u\,.
\end{eqnarray*}

Evaluating at the identity element of $N_4$, $(0,0,0,0)\in\R^4$,
\begin{eqnarray*}
\F_* v &=& v - C_3t\, e - \ve C_2t \, z - \frac{\ve C_2^2 + \bve C_3^2}{2}
   t^2 \, u \, ,\\
\F_* e &=& e + \bve C_3 t\, u\, ,\\
\F_* z &=& z + C_2t \, u \, ,\\
\F_* u &=& u \,.
\end{eqnarray*}
Note that $C_2 = 0$ corresponds to (\ref{d4os}).  On the basis 
$\{v, e, z, u\}$, the matrix of such a $\F_*$ (at the identity) is
$$ \arraycolsep.5em
\left[\begin{array}{cccc}
   1 & 0 & 0 & 0 \\
   -C_3 t & 1 & 0 & 0 \\
   0 & 0 & 1 & 0 \\
   -\half\bve C_3^2 t^2 & \bve C_3 t & 0 & 1 \end{array}\right]. $$
\begin{proposition}
$\F_* \in \Ispl(N_4)$ if and only if $C_2 = C_3 = 0$. $\F_* \in \Iaut(N_4)$ 
if and only if $C_2 = 0$. Thus $\Ispl(N_4) < \Iaut(N_4) < I(N_4)$.\eop
\end{proposition}

Finally, we give the matrix for a $\F_*$ (at the identity) on the basis 
$\{v, e, z, u\}$, as in (\ref{d4os}) but now with $C_3 = 0$ instead. We 
find
$$ \arraycolsep.5em
\left[\begin{array}{cccc}
   1 & 0 & 0 & 0 \\
   0 & 1 & 0 & 0 \\
   -\ve C_2 t & 0 & 1 & 0 \\
   -\half\ve C_2^2 t^2 & 0 & C_2 t & 1 \end{array}\right]. $$
And yet again, these are horocyclic translations as found in \cite{CP1}, 
modulo change of basis.
\end{example}

\section{Lattices and Periodic Geodesics}
\label{latpg}

In this section, we assume that $N$ is rational and let \G\ be a lattice 
in $N$. Since $N$ is 2-step nilpotent, it has nice generating sets of \G; 
{\em cf.}~\cite[(5.3)]{E}. 
\begin{proposition}
If $N$ is a simply connected, 2-step nilpotent Lie group of dimension $n$
with lattice \G\ and with center $Z$ of dimension $m$, then there exists a\/
{\em canonical} generating set $\{\vp_1,\ldots,\vp_n\}$ such that 
$\{\vp_1,\ldots,\vp_m\}$ generate $\G\cap Z$. In particular, $\G\cap Z$ 
is a lattice in $Z$. 
\end{proposition}
From this, formula (\ref{bch}), and \cite[Prop.\,2.17]{R}, one obtains 
as in \cite[(5.3)]{E'} 
\begin{corollary}
Let\label{vl} $N$ be a simply connected, 2-step nilpotent Lie group with 
lattice \G\ and let $\pi:\n\rightarrow\v$ denote the projection. Then
 \begin{enumerate}
  \item $\log\G\cap\z$ is a vector lattice in $\z$;
  \item $\pi(\log\G)$ is a vector lattice in $\v$;
  \item $\G\cap Z=Z(\G)$.\eop
 \end{enumerate}
\end{corollary}
Here, we used the splitting $\n =\U\ds\Z\ds\V\ds\E$ from Section~\ref{defex}
with $\z = \U\ds\Z$ and $\v = \V\ds\E$.

Thus to the compact nilmanifold $\G\bs N$ we may associate two flat 
(possibly degenerate) tori; {\em cf.} \cite[p.\,644]{E'}.
\begin{definition}
With\label{tfb} notation as preceding,
\begin{eqnarray*}
T_{\z} &=& \z/(\log\G\cap\z)\, ,\\
T_{\v} &=& \v/\pi(\log\G)\, .
\end{eqnarray*}
\end{definition}
Observe that $\dim T_{\z} + \dim T_{\v} = \dim\z + \dim\v = \dim\n$.

Next, we apply Theorem 3 from \cite{PS} to our situation. Let $T^m$ denote
the $m$-torus as usual.
\begin{theorem}
Let\label{prin} $m=\dim\z$ and $n=\dim\v$. Then $\G\bs N$ is a principal
$T^m$-bundle over $T^n$.
\end{theorem}
The model fiber $T^m$ can be given a geometric structure from its closed
embedding in $\G\bs N$; we denote this geometric $m$-torus by $T_F$. 
Similarly, we wish to provide the base $n$-torus with a geometric structure 
so that the projection $p_B:\G\bs N\surj T_B$ is the appropriate 
generalization of a pseudoriemannian submersion \cite{O} to (possibly) 
degenerate spaces. Observe that the splitting $\n = \z\ds\v$ induces
splittings $TN = \z N\ds\v N$ and $T(\G\bs N) =
\z(\G\bs N)\ds\v(\G\bs N)$, and that $p_{B*}$ just mods out $\z(\G\bs N)$.
Examining the definition on page 212 of \cite{O}, we see
that the key is to construct the geometry of $T_B$ by defining
\begin{equation}\label{dprs}
p_{B*} : \v_\eta(\G\bs N)\to T_{p_B(\eta)}(T_B) \mbox{ for each }
\eta\in\G\bs N \mbox{ is an isometry }
\end{equation}
and
\begin{equation}\label{dprc}
\nabla^{T_B}_{p_{B*}x}p_{B*}y = p_{B*}\left(\pi\nabla_x y\right)
\mbox{ for all $x,y\in\v=\V\ds\E$, }
\end{equation}
where $\pi:\n\to\v$ is the projection.  Then the rest of the results of
pages 212--213 in \cite{O} will continue to hold, provided that sectional
curvature is replaced by the numerator of the sectional curvature formula
in his 47.\,Theorem at least when elements of $\V$ are involved:
\begin{equation}\label{o47}
\langle R_{T_B}(p_{B*}x,p_{B*}y)p_{B*}y,p_{B*}x\rangle = \langle
R_{\Gamma\bs N}(x,y)y,x\rangle + {\txs\frac{3}{4}}
\langle[x,y],[x,y]\rangle .
\end{equation}
Now $p_B$ will be a pseudoriemannian submersion in the usual sense if and
only if $\U=\V=\{0\}$, as is always the case for Riemannian spaces.

In the Riemannian case, Eberlein showed that $T_F\cong T_{\z}$ and
$T_B\cong T_{\v}$. Observe that it follows from Theorem \ref{kee},
Proposition \ref{oscn}, and (\ref{o47}) that $T_B$ is flat in general
only if $N$ has a nondegenerate center or is flat.
\begin{proposition}
If\label{scTB} $v,v' \in\V$ and $e,e' \in\E$, then
\begin{eqnarray*}
\langle R_{T_B}(v,e)e,v\rangle &=& \quar\langle j(\io v)e,
        j(\io v)e\rangle  ,\\
\langle R_{T_B}(v,v')v',v\rangle &=& \half\langle j(\io v)v',j(\io 
        v')v\rangle - \langle j(\io v)v,j(\io v')v'\rangle\\
 & & {}+\quar\Bigl( \langle j(\io v')v,j(\io v')v\rangle
     +\langle j(\io v)v',j(\io v)v'\rangle\Bigr),\\
K_{T_B}(e,e') &=& 0\,.
\end{eqnarray*}
Here we have suppressed $p_{B*}$ on the left-hand sides for simplicity.\eop
\end{proposition}
Our quaternionic Heisenberg group (Examples \ref{hq} and \ref{hqsc})
provides a simple example of a $T_B$ that is not flat. Indeed, all
sectional curvature numerators vanish except $\langle
R_{T_B}(v_1,v_2)v_2,v_1\rangle = -\bve_1$, where we have again suppressed
$p_{B*}$ for simplicity.
\begin{remark}
Observe\label{sTB} that the torus $T_B$ may be decomposed into a topological
product $T_E \times T_V$ in the obvious way.  It is easy to check that
$T_E$ is flat and isometric to $(\log\G \cap \E)\bs\E$, and that $T_V$ has
a linear connection not coming from a metric and not flat in general.
Moreover, the geometry of the product is ``twisted" in a certain way.  It
would be interesting to determine which tori could appear as such a $T_V$
and how.
\end{remark}

We wish to show that the geometry of the fibers $T_F$ is that of $T_{\z}$.
Thus we now consider the submersion $\G\bs N\surj T_B$ and realize the
model fiber $T^m$ as $(\G\cap Z)\bs Z$ considered as tori only, without
geometry.
\begin{definition}
Let\label{defF} $p_N:N\to\G\bs N$ and $p_Z:Z\to T^m = (\G\cap Z)\bs Z$ be the
natural projections.
Define $F:T^m\to I(\G\bs N)$ by $F\left(p_Z(z)\right)\left(p_N(n)\right) =
p_N(zn)$ for all $z\in Z$ and $n\in N$.
\end{definition}
\begin{proposition}
$F$ is\label{Fiso} a smooth isomorphism of groups $T^m\cong\Ispl_0(\G\bs
N)$, where the subscript 0 denotes the identity component.
\end{proposition}
\begin{proof}
We follow Eberlein \cite[(5.4)]{E'}.  It is easy to check that $F$ is a
well-defined, smooth, injective homomorphism with image in $\Ispl_0(\G\bs
N)$. Thus we need only show that $F$ is surjective.

Let $\psi\in\Ispl_0(\G\bs N)$ and let $\phi_t$ be a path from $1 = \phi_0$
to $\psi=\phi_1$. The covering map $p_N$ has the homotopy lifting
property, so choose a lifting $\tphi_t$ as a path in $\Ispl_0(N)$. Then 
for all $g\in\G$, it follows that $p_N\left(\tphi_t L_g
\tphi_t^{-1}\right) = p_N$ for all $t$. Hence for each $t\in[0,1]$, 
there exists $g_t\in\G$ such that $\tphi_t L_g\tphi_t^{-1} = L_{g_t}$.
Since $g_0 = g$ and $\G$ is a discrete group, it follows that $g_t=g$ for
every $t$, so $L_g$ commutes with $\tphi_t$ for every $g\in\G$ and
$t\in[0,1]$.

From Proposition \ref{sispl}, there exist $n_t\in N$ and $a_t\in O(N)$ such
that $\tphi_t = L_{n_t}a_t$ for all $t$. Now, every $L_g$ commutes with
every $\tphi_t$, so $a_t(g) = n_t^{-1}gn_t$ for all $t$ and $g$. Extension 
from lattices is unique \cite[Thm.\,2.11]{R}, so $a_t =
\Ad{n_t^{-1}}$. By Lemma \ref{iaid}, $a_t$ is the identity and $n_t\in Z$
for all $t$. Thus $\tphi_1 = L_{n_1}$, so from the definition of $\tphi_t$
we obtain $p_N L_{n_1} = p_N\tphi_1 = \phi_1 p_N = \psi p_N$. But this
means $F\left(p_Z(n_1)\right) = \psi$.
\end{proof}
\begin{corollary}
$\Ispl_0(\G\bs N)$ acts\label{fftgo} freely on $\G\bs N$ with complete, flat, 
totally geodesic orbits. 
\end{corollary}
\begin{proof}
By Theorem \ref{prin}, we may identify $\Ispl_0(\G\bs N)$ as the group of
the principal bundle $\G\bs N\surj T_B$, so it acts freely on the total
space.

Since $p_N$ is a local isometry and the $Z$-orbits in $N$ are complete, 
flat, and totally geodesic from Example \ref{ftgZ}, it follows (using the
identification $T^m = (\G\cap Z)\bs Z$ {\em supra}) that the
$\Ispl_0(\G\bs N)$-orbits are complete, flat, and totally geodesic.
\end{proof}
\begin{theorem}
Let\label{TF} $N$ be a simply connected, 2-step nilpotent Lie group with
lattice \G, a left-invariant metric tensor, and tori as in the discussion
following Theorem \ref{prin}.  The fibers $T_F$ of the (generalized)
pseudoriemannian submersion $\G\bs N\surj T_B$ are isometric to $T_{\z}$.
\end{theorem}
\begin{proof}
We follow the proof of Eberlein \cite[(5.5),\,item\,2]{E'}. For each $n\in
N$, define $\psi_n = p_N L_n\exp : \z\to\G\bs N$, and note that it is a
local isometry. Clearly, $\psi_n(z) = \psi_n(z')$ if and only if $z' =
z+\log g$ for some $g\in\G\cap Z$. Hence $\psi_n$ induces an isometric
embedding $\tilde{\psi}_n:T_{\z}\to\G\bs N$. That the image is the
$\Ispl_0(\G\bs N)$-orbit of $p_N(n)$ follows from the proof of Corollary
\ref{fftgo}.
\end{proof}
\begin{corollary}
If\label{fTB} in addition the center $Z$ of $N$ is nondegenerate, then
$T_{B}$ is isometric to $T_{\v}$.\eop
\end{corollary}
The proof is essentially the same as the appropriate parts of the proof of
(5.5) in \cite{E'} and we omit it.

We recall that elements of $N$ can be identified with elements of the
isometry group $I(N)$: namely, $n\in N$ is identified with the isometry
$\phi = L_n$ of left translation by $n$. We shall abbreviate this by
writing $\phi\in N$.
\begin{definition}
We\label{dtrl} say that $\phi\in N$ {\em translates} the geodesic $\g$ by
$\om$ if and only if $\phi\g(t) = \g(t+\om)$ for all $t$. If $\g$ is a
unit-speed geodesic, we say that $\om$ is a\/ {\em period} of $\phi$.
\end{definition} 
Recall that unit speed means that $|\dot{\g}| = \left|\langle\dot{\g},
\dot{\g}\rangle\right|^{\frac{1}{2}} = 1$.
Since there is no natural normalization for null geodesics, we do not
define periods for them. In the Riemannian case and in the timelike
Lorentzian case in strongly causal spacetimes \cite{BE}, unit-speed
geodesics are parameterized by arclength and this period is a translation
distance. If $\phi$ belongs to a lattice $\G$, it is the length of a
closed geodesic in $\G\bs N$.
\begin{remark}
Note\label{rcc} that it follows from Corollary \ref{pcc} that if $\phi =
\exp(a^* + x^*)$ translates a geodesic $\g$ with $\g(0) = 1 \in N$, then
$a^* + x^*$ and $\dg(0)$ are of the same causal character.
\end{remark}

In general, recall that if $\g$ is a geodesic in $N$ and if $p_N:N\surj
\G\bs N$ denotes the natural projection, then $p_N\g$ is a periodic
geodesic in $\G\bs N$ if and only if some $\phi\in\G$ translates $\g$.
We say {\em periodic\/} rather than {\em closed\/} here because in
pseudoriemannian spaces it is possible for a null geodesic to be closed
but not periodic. If the space is geodesically complete or Riemannian,
however, then this does not occur ({\em cf.}~\cite{O}, p.\,193); the
former is in fact the case for our 2-step nilpotent Lie groups by
Theorem \ref{gt}.  Further recall that free homotopy classes of closed
curves in $\G\bs N$ correspond bijectively with conjugacy classes in $\G$.
\begin{definition}
Let\label{pC} $\sC$ denote either a nontrivial, free homotopy class of
closed curves in\/ $\G\bs N$ or the corresponding conjugacy class in\/
$\G$. We define $\wp(\sC)$ to be the set of all periods of periodic
unit-speed geodesics that belong to $\sC$.
\end{definition}
In the Riemannian case, this is the set of lengths of closed geodesics in
$\sC$, frequently denoted by $\ell(\sC)$.
\begin{definition}
The\label{psp} {\em period spectrum\/} of\/ $\G\bs N$ is the set
$$ \specp(\G\bs N) = \bigcup_{\sC}\wp(\sC)\,,$$
where the union is taken over all nontrivial, free homotopy
classes of closed curves in $\G\bs N$.
\end{definition}
In the Riemannian case, this is the length spectrum $\specl(\G\bs N)$.
\begin{example}
Similar\label{sft} to the Riemannian case, we can compute the period
spectrum of a flat torus $\G\bs\R^m$, where $\G$ is a lattice (of maximal
rank, isomorphic to ${\Bbb Z}^m$).  Using calculations related to those of
\cite[pp.\,146--8]{BGM} in an analogous way as for finding the length
spectrum of a Riemannian flat torus, we easily obtain
$$ \specp(\G\bs\R^m) = \{ |g|\ne 0 \mid g\in\G\}\,. $$
It is also easy to see that the nonzero d'Alembertian spectrum is related
to the analogous set produced from the dual lattice $\G^*$ as multiples by
$\pm 4\pi^2$, almost as in the Riemannian case.
\end{example}

As in this example, simple determinacy of periods of unit-speed geodesics
helps make calculation of the period spectrum possible purely in terms of
$\log\G \subseteq \n$.  (See Theorem \ref{upd} for another example.)  Thus
we begin with the following observation.
\begin{proposition}
Let\label{vup} $\phi = \exp(a^* + v^* + e^*)$ translate the unit-speed
geodesic $\g$ by $\om > 0$.  If $v^* \ne 0$, then the period $\om$ is
simply determined.
\end{proposition}
\begin{proof}
We may assume $\g(0) = 1 \in N$ and $\dg(0) = a_0 + v_0 + e_0$.
From Theorem \ref{ige}, $v^* = v(\om) = \om\,v_0$.
\end{proof}
In attempting to calculate the period spectrum then, we can focus our
attention on those cases where the $\V$-component is zero.

From now on, we assume that $N$ is a simply connected, 2-step nilpotent
Lie group with left-invariant pseudoriemannian metric tensor
$\langle\,,\rangle$.  Note that non-null geodesics may be taken to be of
unit speed.  Most nonidentity elements of $N$ translate some geodesic, but
not necessarily one of unit speed; {\em cf.} \cite[(4.2)]{E'}.
\begin{proposition}
Let\label{tsg} $N$ be a simply connected, 2-step nilpotent Lie group with
left-invariant metric tensor $\langle\,,\rangle$ and $\phi\in N$ not the
identity. Write $\log\phi = a^* +x^*\in\z\ds\v$ and assume that $x^* \perp
[x^*,\n]$.  Let $a'$ be the component of $a^*$ orthogonal to $[x^*,\n]$ in
$\z$ and choose $\xi\in\n$ such that $a'=a^* + [x^*,\xi]$.  Set $\om^*
= |a' +x^*|$ if $a'+x^*$ is not null and set $\om^* = 1$ otherwise.  Then
$\phi$ translates the geodesic 
$$ \g(t)= \exp(\xi)\,\exp\left(\frac{t}{\om^*}(a'+x^*)\right) $$
by $\om^*$, and $\g$ is of unit speed if $a'+x^*$ is not null.
\end{proposition}
\begin{proof}
Let $n=\exp(\xi)$ and set $\phi^*=n^{-1}\phi n=\exp(a'+x^*)$ and $\g^*(t)
= n^{-1}\g(t)$. Then $\phi\g(t) = \g(t+\om^*)$ is equivalent to
$\phi^*\g^*(t) = \g^*(t+\om^*)$, and the latter is routine to verify using
(\ref{bch}).

Now $\g$ is a geodesic if and only if $\g^*$ is, and it is easy to check
directly that $\del{a'+x^*}(a'+x^*) = 0$ is equivalent to $\langle
a'+x^*, [x^*,\n]\rangle = 0$.
\end{proof}
Note that the $\U$ components of $a^*$ and $a'$ in fact coincide.  Also
note that if we further decompose $x^* = v^* + e^*$, the result applies to
every nonidentity element $\phi$ with $v^* = 0$.  In particular, when the
center is nondegenerate this is every nonidentity element.
\begin{corollary}
When $\n$ is nonsingular and $\phi\notin Z$, we may take $a'=0$ in
Proposition \ref{tsg}.
\end{corollary}
\begin{proof}
Because then $a^*\in[x^*,\n]$ and $x^*\ne 0$.
\end{proof}

Now we give some general criteria for an element $\phi$ to translate a
geodesic $\g$; {\em cf.} \cite[(4.3)]{E'}.  We use a $J$ as in the passage
following Definition \ref{defJ}, and $x_1$, $x_2$, $y_1$, and $y_2$ as
given just before Theorem \ref{ige}.
\begin{proposition}
Let\label{gctg} $\phi\in N$ and write $\phi = \exp(a^* + x^*)$ for suitable
elements $a^*\in\z$ and $x^*\in\v$. Let $\g$ be a geodesic with $\g(0) =
n\in N$ and $\g(\om) = \phi n$, and let $\dg(0) = L_{n*}(a_0 + x_0)$ for
suitable elements $a_0\in\z$ and $x_0 = v_0+e_0\in\v$.  Let $n^{-1}\g(t)=
\exp\left(a(t) + x(t)\right)$ where $a(t)\in\z$ and $x(t)\in\v$ for all
$t\in\R$ and $a(0) = x(0) = 0$. Then the following are equivalent:
\begin{enumerate}
\item $\left.\begin{array}{rcl}
            x(t+\om) &=& x(t) + x^* \\
             a(t+\om) &=& a(t) + a^* + \frac{1}{2}\left[x^*,x(t)\right]
             \end{array}\right\}$ for all $t\in\R$ and some $\om>0$;

\item $\g(t+\om) = \phi\g(t)$ for all $t\in\R$ and some $\om>0$;

\item $e^{\om J}$ fixes $e_1 + y_1 + x_2 = e_0 + y_1 + J^{-1}y_2$.
\end{enumerate}
\end{proposition}
\begin{proof}
As before, we may assume without loss of generality that $n = 1\in N$.
Items 1 and 2 are equivalent by formula (\ref{bch}).  The following lemma
shows that item 1 implies item 3.
\begin{lemma}
As\label{lem} in the preamble to Theorem \ref{ige}, write $\E$ as an
orthogonal direct sum $\E_1\ds\E_2$, where $\E_1 = \ker J$, and use $x_1$,
$x_2$, $y_1$, and $y_2$ as given there. Then $x^* = \om x_1 + \half\om^2 
y_1$ and $e^{\om J}$ fixes $x_2$.
\end{lemma}
\begin{proof}
By Theorem \ref{ige}, we have $x(k\om) = k\om x_1 + \left(e^{k\om J} -
I\right) J^{-1}x_2 + \half k^2\om^2 y_1$ for every positive integer $k$.
By induction, from item 1 in the statement of the proposition we obtain
$x(k\om) = k x^*$ for every $k$. Decomposing $x^* = v^* + e^*_1 + e^*_2
\in \V\ds\E_1\ds\E_2$, this yields $v^* = \om v_0$, $e^*_1 = \om e_1 +
\half\om^2 y_1$, and
\begin{equation}
k( e^*_2 + \om J^{-1}y_2) = \left(e^{k\om J} - I\right)J^{-1}x_2 
\label{unbdd}
\end{equation}
for every $k$.

Now, $e^{\om J}$ is an element of the identity component of the
pseudorthogonal group of isometries of $\langle\,,\rangle$, and as such
can be decomposed into a product of reflections, ordinary rotations, and
boosts. With respect to appropriate coordinates, which may be different
from our standard choice, a boost will have a matrix of the form
$$ \left[\begin{array}{cc}
         \cosh s & \sinh s\\
         \sinh s & \cosh s \end{array}\right] $$
on some pair of basis vectors, for some $s\in\R$.

If $e^{\om J}$ is composed only of reflections and ordinary rotations,
then the right-hand side of (\ref{unbdd}) is uniformly bounded in $k$
(say, with respect to the positive definite $\langle\,,\io\rangle$) while
the left-hand side is unbounded, so $e^*_2 + \om J^{-1}y_2 = 0$.  On the
other hand, if $e^{\om J}$ is a pure boost, then the right-hand side grows
exponentially in $k$ while the left-hand side grows but linearly, and
again we obtain $e^*_2 + \om J^{-1}y_2 = 0$.  The Lemma now follows from
this and (\ref{unbdd}) for $k=1$.
\end{proof}
The proof that item 3 implies item 2 is the same as the relevant part of
the proof of (4.3) in \cite{E'}.
\end{proof}
\begin{corollary}
When\label{gctgc} in addition $v^* = v_0 = 0$, the following are also
equivalent to the three items in Proposition \ref{gctg}.
\begin{enumerate}
\item $\dg(0)$ is orthogonal to the orbit $Z_{e^*}n$, where 
$Z_{e^*} = \exp\left([e^*,\n]\right) \subseteq Z$;

\item $\dg(\om)$ is orthogonal to the orbit $Z_{e^*}\phi n$.
\end{enumerate}
\end{corollary}
\begin{proof}
Note that under this hypothesis, $x^* = e^*$.  Now Lemma \ref{lem}
implies that $J(e^*) = 0$, and this is now equivalent to $z_0 \perp
[e^*,\n]$.  Thus the relevant parts of the proof of (4.3) in \cite{E'}
apply {\em mutatis mutandis.}
\end{proof}

We also obtain the following results as in Eberlein
\cite[(4.4),\,(4.9)]{E'}.  Note that we assume that $v^* = v_0 = 0$ in the
first, but that this is automatic in the second.
\begin{corollary}
Let\label{tsg1} $\phi\in N$ and write $\phi = \exp(a^* + e^*)$ for unique
elements $a^*\in\z$ and $e^*\in\E$. Let $n\in N$ be given and write $n =
\exp(\xi)$ for a unique $\xi\in\n$. Then the following are equivalent:
\begin{enumerate}
\item There exists a geodesic $\g$ in $N$ with $\g(0) = n$ such
that $\phi\g(t) = \g(t + \om)$ for all $t\in\R$ and some $\om > 0$.

\item There exists a geodesic $\g^*$ in $N$ with $\g^*(0) = 1$,
$\dg^*(0)$ is orthogonal to $[e^*,\n]$, and $\g^*(\om)
= \exp\left([e^*, \xi]\right)\phi$ for some $\om > 0$.\eop
\end{enumerate}
\end{corollary}
\begin{corollary}
Let\label{tsg2} $1\ne\phi\in Z$ and $\g$ be any geodesic such that $\g(\om)
= \phi\g(0)$ for some $\om > 0$.  Then $\phi\g(t) = \g(t + \om)$ for all
$t\in\R$.\eop
\end{corollary}

In the flat 2-step nilmanifolds of Theorem \ref{e0f}, we can calculate the
period spectrum completely.
\begin{theorem}
If\label{upd} $[\n,\n] \subseteq \U$ and $\E=\{0\}$, then $\specp(M)$ can be
completely calculated from $\log\G$ for any $M = \G\bs N$.
\end{theorem}
\begin{proof}
Let $\phi$ translate a unit-speed geodesic $\g$ by $\om > 0$.  As usual,
we may as well assume that $\g(0) = 1 \in N$.  Write $\log\phi = a^* +
v^*$ and $\dg(0) = a_0 + v_0$.  From Corollary \ref{iged} we get $v^* =
v(\om) = \om\,v_0$, $z^* = z(\om) = \om\,z_0$, and $u^* = u(\om) =
\om\,u_0 + \half\om^2\,\sJ v_0$.  Note that $\om^2\,\sJ v_0 =
\om^2\,\add{v_0}{(z_0 + v_0)} = \om^2\,\add{v_0}{v_0} = \add{v^*}{v^*}$.
Substituting and rearranging, we obtain
$$ \ve\om^2 = 2\langle u^*,v^*\rangle + \langle z^*,z^*\rangle , $$
where $\pm1 = \ve = \langle\dg(0),\dg(0)\rangle = 2\langle u_0,v_0\rangle
+ \langle z_0,z_0\rangle$.
\end{proof}
Thus we see again, as mentioned after Corollary \ref{iged}, just how much
these flat, 2-step nilmanifolds are like tori.  All periods can be
calculated purely from $\log\G \subseteq \n$, although some will not show
up from the tori in the fibration.
\begin{corollary}
$\specp(T_B)$ (respectively, $T_F$) is\label{psd} $\cup_{\sC}\,\wp(\sC)$
where the union is taken over all those free homotopy classes\/ $\sC$ of
closed curves in $M = \G\bs N$ that\/ {\em do not} (respectively,\/ {\em
do}) contain an element in the center of\/ $\G \cong \pi_1(M)$, except for
those periods arising only from unit-speed geodesics in $M$ that project to
null geodesics in both $T_B$ and $T_F$.\eop
\end{corollary}
We note that one might consider using this to assign periods to some null
geodesics in the tori $T_B$ and $T_F$.

When the center is nondegenerate, we obtain results similar to Eberlein's
\cite[(4.5)]{E'}.
\begin{proposition}
Assume\label{e4.5} $\U = \{0\}$. Let $\phi \in N$ and write $\log\phi =
z^* + e^*$.  Assume $\phi$ translates the unit-speed geodesic $\g$ by $\om
> 0$.  Let $z'$ denote the component of $z^*$ orthogonal to $[e^*,\n]$.
Let $n = \g(0)$ and set $\om^* = |z' + e^*|$. Let $\dg(0) = L_{n*}(z_0 +
e_0)$ and use $J$, $e_1$, and $e_2$ as in Corollary \ref{igend} (see also
just before Theorem \ref{ige}).  Then
\begin{enumerate}
\item $|e^*| \le \om$. In addition, $\om < \om^*$ for timelike
   (spacelike) geodesics with $\om z_0 - z'$ timelike (spacelike), and
   $\om > \om^*$ for timelike (spacelike) geodesics with $\om z_0 - z'$
   spacelike (timelike).

\item $\om = |e^*|$ if and only if $\g(t) = \exp\left(t\,e^*\!/|e^*|
   \right)$ for all $t \in \R$.

\item $\om = \om^*$ if and only if $\om z_0 - z'$ is null.  If moreover
   $\om^* z_0 = z'$, then $e_2 = 0$ if and only if
\begin{enumerate}
   \item $\g(t) = n\,\exp\left( t\,(z' + e^*)/\om^*\right)$ for all
      $t \in \R$.

   \item $z' = z^* + [e^*,\xi]$ where $\xi = \log n$.
\end{enumerate}
\end{enumerate}
\end{proposition}
Although $\om^*$ need not be an upper bound for periods as in the
Riemannian case, it nonetheless plays a special role among all periods,
as seen in item 3 above, and we shall refer to it as the {\em
distinguished\/} period associated with $\phi \in N$.  When the center is
definite, for example, we do have $\om \le \om^*$.
\begin{proof}
As usual, we may assume that $\g(0) = 1 \in N$.  Note that this replaces
$\phi$ as given in the statement with $n^{-1}\phi n$ and $\g$ with
$n^{-1}\g$.

For the first part of item 1, since $|\dg(0)| = 1$ there exists an
orthonormal basis of $\n$ having $\dg(0)$ as a member.  (This may well be
a different basis from our usual one.)  Fix one such basis, and consider
the positive-definite inner product with matrix $I$ on this basis.  Let
$\|\cdot\|$ denote the norm associated to this positive-definite inner
product. By Lemma \ref{lem}, $e^* = \om e_1$.  Then $|e_1| \le \|e_1\| \le
\|\dg(0)\| = 1$ so $|e^*| =
\om |e_1| \le \om$.

For the rest of item 1, we begin with Corollary \ref{igend} and get
\begin{eqnarray*}
e(t) &=& t\,e_1 + \left(e^{tJ} - I\right)J^{-2}e_2\,,\\
z(t) &=& t\left( z_0 + \half\left[ e_1, \left(e^{tJ}+I\right)J^{-1}e_2 
         \right] \right) + z_2(t) + z_3(t)\,.
\end{eqnarray*}
By Lemma \ref{lem}, $e^* = \om e_1$ and $e^{\om J}e_2 = e_2$. Inspecting
the formula for $z_2(t)$ in Corollary \ref{igend}, we find $z_2(\om) = 0$.
Thus
\begin{eqnarray*}
z^* = z(\om) &=& \om\left(z_0 + [e_1,J^{-1}e_2]\right) + z_3(\om)\\
   &=& \om\,z_0 + [e^*,J^{-1}e_2] + \half\int_0^{\om} \left[ e^{sJ}
J^{-1}e_2, e^{sJ}e_2\right] ds\,.
\end{eqnarray*}
By item 1 of Corollary \ref{gctgc}, $z_0 \perp [e^*,\n]$.  Then
$$ \langle z',z_0\rangle = \langle z^*,z_0\rangle = \om\langle
z_0,z_0\rangle + \half\int_0^{\om}\left\langle\left[e^{sJ}J^{-1}e_2,
e^{sJ}e_2\right], z_0\right\rangle ds\,.$$
Recall that $J$ is skewadjoint with respect to $\langle\,,\rangle$ (whence
$e^{tJ}$ is an isometry of $\langle\,,\rangle$ for all $t$), that $J$
commutes with every $e^{tJ}$ (whence so does $J^{-1}$), and that $Jx =
\add{x}{z_0}$. We compute
\begin{eqnarray*}
\left\langle\left[e^{sJ}J^{-1}e_2, e^{sJ}e_2\right], z_0\right\rangle
&=& -\left\langle\left[e^{sJ}e_2, e^{sJ}J^{-1}e_2\right], z_0\right\rangle\\
&=& -\left\langle J^{-1}e^{sJ}e_2, Je^{sJ}e_2\right\rangle\\
&=& \left\langle e^{sJ}e_2, e^{sJ}e_2\right\rangle\\
&=& \langle e_2, e_2\rangle.
\end{eqnarray*}
Therefore,
\begin{equation}
\langle z', z_0\rangle = \om\langle z_0, z_0\rangle +
\frac{\om}{2}\langle e_2, e_2\rangle .\label{one}
\end{equation}

Now $|\dg(0)| = 1$ so $\ve = \langle z_0,z_0\rangle + \langle e_1,
e_1\rangle + \langle e_2,e_2\rangle$, where $\ve = \pm 1$ as usual.
Substituting in (\ref{one}) for $\langle e_2,e_2\rangle$, we obtain
$$ \langle z',z_0\rangle = \frac{\om}{2}\bigl( \ve + \langle
   z_0,z_0\rangle\bigr) - \frac{\om}{2}\frac{\langle
   e^*,e^*\rangle}{\om^2} $$
so
$$ \langle e^*,e^*\rangle - \ve\om^2 = \om^2\langle z_0,z_0\rangle
-2\om\langle z',z_0\rangle. $$
Adding $\langle z',z'\rangle$ to both sides, we get
\begin{equation}
\langle z'+e^*, z'+e^*\rangle - \ve\om^2 = 
   \langle \om z_0-z', \om z_0-z'\rangle.\label{two}
\end{equation}
There are several cases: $\ve$ is 1 or $-1$ and $\om z_0-z'$ is timelike,
spacelike, or null.  If $\om z_0-z'$ is null, then $|\langle z'+e^*,
z'+e^*\rangle| = \om^2 > 0$ and $\om = |z'+e^*|$.  If $\ve = 1$ and $\om
z_0-z'$ is timelike, or if $\ve = -1$ and $\om z_0 - z'$ is spacelike,
then $\ve\langle z'+e^*, z'+e^*\rangle > \om^2 > 0$ whence $\om <
|z'+e^*|$.  If $\ve = 1$ and $\om z_0-z'$ is spacelike, or if $\ve = -1$
and $\om z_0 - z'$ is timelike, then it follows similarly that $\om > |z'
+ e^*|$.  This completes the proof of item 1.

Now we prove item 2.  If $\g$ is as given there, then $\exp(z^* + e^*) =
\phi = \g(\om) = \exp(\om\,e^*\!/|e^*|)$ whence $\om = |e^*|$.
Conversely, assume $\om = |e^*|$ and consider the associated
positive-definite inner product $\langle\cdot,\io\cdot\rangle$.  Changing
the basis of $\E$ if necessary, we may assume that $\Z$, $\E_1$, and
$\E_2$ are mutually orthogonal with respect to both $\langle\,,\rangle$
and $\langle\cdot,\io\cdot\rangle$.  Let $\|\cdot\|$ now denote the norm
for $\langle\cdot,\io\cdot\rangle$.  Then
\begin{equation}
\|\dg(0)\|^2 = \|z_0\|^2 + \|e_1\|^2 + \|e_2\|^2\label{three}
\end{equation}
so $\|\dg(0)\|^2 = \|e_1\|^2$ if and only if $\dg(0) = e_1 = e^*\!/|e^*|$.
But now $\g$ has the same initial data as $\exp(t\,e^*\!/|e^*|)$, so by
uniqueness they must coincide.

Finally, we prove the last part of item 3; the first part is immediate
from the last part of the proof of item 1 above.  So assume $\om^*z_0 - z'
= 0$ or $z_0 = z'\!/\om^*$. Continue with the immediately previous
positive-definite norm $\|\cdot\|$ and basis of $\E$.  Substituting in
(\ref{three}) we get
\begin{eqnarray*}
\|\dg(0)\|^2 &=& \frac{\|z'\|^2}{(\om^*)^2} + \frac{\|e^*\|^2}{(\om^*)^2}
+ \|e_2\|^2\\
   &=& \frac{\|z'+e^*\|^2}{(\om^*)^2} + \|e_2\|^2
\end{eqnarray*}
whence
$$\dg(0) = \frac{z'+e^*}{\om^*}$$
if and only if $e_2 = 0$.
\end{proof}
\begin{corollary}
Assume\label{e4.6} the center is nondegenerate.  Let $\phi \in N$ with
$\phi \notin Z$ and suppose that $z^* \in [e^*,\n]$. Then
\begin{enumerate}
\item If $\phi$ translates a timelike (spacelike) geodesic with $z_0$
nonspacelike (nontimelike), then $\phi$ has the unique period $|e^*|$.

\item Let $\g$ be a unit-speed geodesic in $N$ with $\g(0) = n =
\exp(\xi)$ for a unique $\xi \in \n$. Then $\phi$ translates $\g$ by the
unique period $|e^*| > 0$ if and only if $[\xi,e^*] = z^*$ and $\g(t) =
n\,\exp(t\,e^*\!/|e^*|)$ for all $t \in \R$.
\end{enumerate}
In particular, this applies to all noncentral\/ $\phi \in N$ if\/ $\n$ is
nonsingular.\eop
\end{corollary}
The proof follows that of \cite[(4.6)]{E'} {\em mutatis mutandis\/} and we
omit the details. We do note that from item 2 of Proposition \ref{e4.5},
using Lemma \ref{dexp}, we obtain $\exp(e^* + \half[\xi,e^*]) =
n\,\exp(e^*) = \g(|e^*|) = \phi\,n = \exp(z^* + e^*)\exp(\xi) = \exp(z^* +
e^* + \half[e^*,\xi])$, thus avoiding the use of item 3 here.  Anent the
last comment, note that if $\n$ is nonsingular then in fact $z_0 = 0$ in
item 1, because $z_0 \perp [e^*,\n] = \z$.

In view of the comment following Proposition \ref{e4.5}
and Corollary \ref{e4.6}, the following definitions make sense at
least for $N$ with a nondegenerate center.
\begin{definition}
Let\label{MpC} $\sC$ denote either a nontrivial, free homotopy class of
closed curves in\/ $\G\bs N$ or the corresponding conjugacy class in\/
$\G$. We define $\wp^*(\sC)$ to be the distinguished periods of periodic
unit-speed geodesics that belong to $\sC$.
\end{definition}
\begin{definition}
The\label{Dpsp}\/ {\em distinguished period spectrum} of\/ $\G\bs N$ is the
set
$$ \Dspecp(\G\bs N) = \bigcup_{\sC}\wp^*(\sC)\,,$$
where the union is taken over all nontrivial, free homotopy classes of
closed curves in $\G\bs N$.
\end{definition}
Then as an immediate consequence of the preceding corollary, we get this
result.
\begin{corollary}
Assume\label{psnd} the center is nondegenerate.  If $\n$ is nonsingular,
then $\specp(T_B)$ (respectively, $T_F$) is precisely the period spectrum
(respectively, the distinguished period spectrum) of those free homotopy
classes\/ $\sC$ of closed curves in $M = \G\bs N$ that\/ {\em do not}
(respectively,\/ {\em do}) contain an element in the center of\/ $\G \cong
\pi_1(M)$, except for those periods arising only from unit-speed geodesics
in $M$ that project to null geodesics in both $T_B$ and $T_F$.\eop
\end{corollary}

\section{Constructing Lie algebras of \ph-type}
\label{cph}

Let us suppose that $\n = \z\ds\v = \U\ds\Z\ds\V\ds\E$ is a vector space
(not a Lie algebra), for which we have a basis consisting of
$$ \begin{array}{lllllll}
\{u_i\} & \mbox{a basis for} & \U\, , & \qquad &
\{z_\alpha\} & \mbox{a basis for} & \Z\,, \\
\{v_i\} & \mbox{a basis for} & \V\, , & \qquad  &
\{e_a\} & \mbox{a basis for} & \E\,.
\end{array} $$
Assume that there exists a nontrivial inner product $\langle\,,\rangle$ on 
$\n$ which is given with respect this basis by
$$
\langle u_i,v_i \rangle = 1\,, \quad \langle z_\alpha ,z_\alpha \rangle =
\ve_\alpha\,, \quad \langle e_a,e_a \rangle = \bar{\ve}_a\,.
$$
Here, each $\ve$-symbol is $\pm1$ independently, making the bases for $\Z$
and $\E$ orthonormal.

Then we  define an involution $\io \colon \n \to \n$ with respect to
this basis by setting:
$$
\io (u_i)=v_i\,, \quad \io (v_i)=u_i\,, \quad \io (z_\alpha )=\ve_\alpha \,
z_\alpha\, , \quad \io (e_a )= \bar{\ve}_a \, e_a\,.
$$
Thus $\io (\U )=\V$, $\io (\V )=\U$, $\io (\Z )=\Z$, $\io (\E )=\E$, and
$\io^2 = I$. Moreover, $\io$ is self-adjoint with respect to the inner 
product: 
$\langle \io x,y \rangle = \langle x, \io y \rangle$ for every $x,y \in \n$.
Therefore $\io$ is an isometry of \n.
Note that $\langle x, \io x \rangle = 0$ if and only if $x = 0$.

Let us now also assume that there is given a linear mapping
$$
j:\U\ds\Z\rightarrow\End(\V\ds\E)
$$
satisfying the conditions
\begin{eqnarray}
\langle j(a)x,\io j(a)x \rangle &=& \langle a,\io a \rangle
          \langle x,\io x \rangle \label{id-1}\\
j(a)^2 &=& -\langle a,\io a \rangle I \label{id-2}
\end{eqnarray}
for every $a\in\U\ds\Z$ and $x\in\V\ds\E$.
 
Then there exists on $\n$ a canonical Lie algebra structure of \ph-type,
obtained as follows.
 
First, we note that polarizing (\ref{id-1}) and (\ref{id-2}) yields the
relations
\begin{eqnarray}
\langle j(a)x,\io j(b)x \rangle &=& \langle a,\io b \rangle
          \langle x,\io x \rangle\\
\langle j(a)x,\io j(a)y \rangle &=& \langle a,\io a \rangle
          \langle x,\io y \rangle \label{id-3}
\end{eqnarray}
for every $a,b\in\U\ds\Z$ and $x,y\in\V\ds\E$. Then
it is easy to prove that the map $j$ is $\io$-skewsymmetric:
$$
\langle j(a)x,\io y \rangle + \langle \io x,j(a)y \rangle = 0
$$
for every $a\in\U\ds\Z$ and $x,y\in\V\ds\E$. Indeed, if $a=0$ this is
trivial, and if $a\ne 0$ then from (\ref{id-1}),
(\ref{id-2}), and (\ref{id-3}), we get
$$
\begin{array}{rcl}
\langle j(a)x,\io y \rangle &=& \langle j(a)x,\io j(a)^2 \left(
  -\langle a,\io a \rangle^{-1}y \right) \rangle \\
&=& -\langle a,\io a \rangle^{-1} \langle j(a)x,\io j(a)\left(
     j(a)y \right) \rangle \\
&=& -\langle a,\io a \rangle^{-1} \langle a,\io a \rangle
     \langle x,\io j(a)y \rangle \\
&=& -\langle \io x,j(a)y \rangle .
\end{array}
$$
 
Now we define a bilinear map
$$
[\, ,\, ] : (\V\ds\E )\ds(\V\ds\E ) \longrightarrow \U\ds\Z
$$
as follows:
$$
\langle [x,y],\io a \rangle = \langle j(a)x, \io y \rangle
$$
for every $x,y\in\V\ds\E$ and $a\in\U\ds\Z$. 
We remark that $[x,y]\in\U\ds\Z$ is well defined.
\begin{proposition}
This map $[\, ,\, ]$ is skewsymmetric: $[x,y]=-[y,x]$ for every 
$x,y\in\V\ds\E$.
\end{proposition}
\begin{proof} We compute:
\begin{eqnarray*}
\langle [x,y],\io a \rangle &=& \langle j(a) x,\io y \rangle\\
&=& -\langle \io x,j( a) y \rangle\\
&=& -\langle [y,x],\io a \rangle . \eop
\end{eqnarray*}
\end{proof}
 
We can now define a Lie algebra structure on $\n$ by extending this map
to all of $\n$ {\em via}
$$
[a+x,b+y]=[x,y]\in\U\ds\Z
$$
for every $a,b\in\U\ds\Z$, and $x,y\in\V\ds\E$. Clearly, we obtain
\begin{theorem}
Endowed with this structure, $\n$ is a $2$-step nilpotent
Lie algebra of \ph-type with center $\U\ds\Z$.\eop
\end{theorem}
 
We illustrate this construction with two simple examples.
\begin{example}
Assume $\U =\V=\{0\}$, $\dim\Z =1$, and $\dim\E=2$, so that $\dim\n =3$.
Choose bases $\{z\}$ of $\Z$ and $\{e_1,e_2\}$ of $\E$, and
nontrivial inner products
$$ \ve=\langle z,z\rangle ,\quad\bve_i=\langle e_i,e_i\rangle . $$
Suppose that the mapping
$$ j:\Z\rightarrow\End(\E) $$
is given with respect to the basis $\{e_1,e_2\}$ as
$$ j(z) = \left[
\begin{array}{cc}
0 & -1\\
1 & 0
\end{array}
\right]. $$
Then
$$ \langle [e_1,e_2],\io z \rangle = \ve \langle [e_1,e_2],z \rangle , $$
and, on the other hand,
$$ \langle j(z ) e_1,\io e_2 \rangle = \bve_2 \langle j(z)e_1,e_2
\rangle = \bve_2 \langle  e_2,e_2 \rangle = 1\,. $$
Hence, $\ve\langle [e_1,e_2],z \rangle =1$ so $[e_1,e_2]=z$.

Therefore, $\n$ is the $3$-dimensional Lie algebra with structure equation
$[e_1,e_2]=z$; that is, the $3$-dimensional Heisenberg algebra with non-null
center.
\end{example}
\begin{example}
Assume $\dim\U =\dim\V=1$, $\Z =\{0\}$ and $\dim\E=1$, so that $\dim\n =3$
again. Choose bases $\{u\}$ for $\U$, $\{v\}$ for $\V$, and $\{e\}$ for $\E$,
and nontrivial inner products
$$ \langle u,v \rangle = 1\, , \quad \langle e,e \rangle = \bve\, . $$
Suppose that the mapping
$$ j:\U\rightarrow\End(\V\ds\E) $$
is given with respect to the basis $\{v,e\}$ as
$$ j(u)=\left[
\begin{array}{cc}
0 & -1 \\
1 & 0
\end{array}
\right]. $$
Then
$$ \langle [v,e],\io u \rangle = \langle [v,e],v \rangle , $$
and, on the other hand,
$$ \langle j(u ) e,\io v \rangle =  \langle j(u)e,u
\rangle =  \langle v,u \rangle = 1\, . $$
Hence, $\langle [v,e],v \rangle =1$ so $[v,e]=u$.
 
Therefore, $\n$ is the $3$-dimensional Lie algebra with structure equation
$[v,e]=u$; that is, the $3$-dimensional Heisenberg algebra with null
center.
\end{example}

%\newpage
\section*{Acknowledgments}
Once again, Parker wishes to thank the Departamento at Santiago for its
fine hospitality.  He also thanks WSU for a Summer Research Fellowship in
1994 during which part of this work was done, and for a Sabbatical Leave
in 1998 during which it was finished.

\end{document}